\documentstyle[amsfonts, amssymb, latexsym, epsfig, epic, eepic, 12pt]{amsart}

\newtheorem{theorem}{Theorem}[section]
\newtheorem{lemma}[theorem]{Lemma}
\newtheorem{proposition}[theorem]{Proposition}
\newtheorem{corollary}[theorem]{Corollary}
\newtheorem{definition}[theorem]{Definition}

\newtheorem{question}[theorem]{Question}

\def\Z{{\mbox{\rm\kern.25em
\vrule width.03em height0.57ex depth0ex
\kern.033em
\vrule width.03em height1.52ex depth-0.96ex \kern-.338em Z}}}
\def\R{{\mbox{\rm I\kern-.22em R}}}
\def\N{{\mbox{\rm I\kern-.22em N}}}

\def\D{{\bf D}}
\def\T{{\bf T}}
\def\supp{{\rm supp}}
\def\size{{\rm size}}

\def\energy{{\rm energy}}
\def\max{{\rm max}}

\def\M{{\cal{M}}}

\def\k{{\kappa}}

\def\L{{\cal{L}}}
\def\D{{\cal{D}}}
\def\J{{\cal{J}}}
\def\I{{\cal{I}}}

\def\dist{{\rm dist}}

\def\111{\gamma}

\def\be#1{\begin{equation}\label{#1}}
\def\bas{\begin{align*}}
\def\eas{\end{align*}}
\def\bi{\begin{itemize}}
\def\ei{\end{itemize}}
\newenvironment{proof}{\noindent {\bf Proof} }{\endprf\par}
\def \endprf{\hfill  {\vrule height6pt width6pt depth0pt}\medskip}
\def\emph#1{{\it #1}}

\title{Paraproducts with flag singularities I. A case study}

\author{Camil Muscalu}
\address{Department of Mathematics, Cornell University, Ithaca, NY 14853}
\email{camil@@math.cornell.edu}

\include{psfig}
\begin{document}

\begin{abstract}
In this paper we prove $L^p$ estimates for a tri-linear operator, whose symbol
is given by the product of two standard symbols, satisfying the well known 
Marcinkiewicz-H\"{o}rmander-Mihlin condition. Our main result contains in 
particular the classical Coifman-Meyer theorem. This tri-linear operator is 
the simplest example of a large class of multi-linear operators,
which we called {\it paraproducts with flag singularities}.

\end{abstract}

\maketitle

\section{Introduction}

The purpose of the present article is to start a systematic study of the
$L^p$ boundedness properties of a new class of multi-linear operators which 
we named {\it paraproducts with flag singularities}.

For any $d\geq 1$ let us denote by $\cal{M}$$(\R^d)$ the set of all bounded 
symbols $m\in L^{\infty}(\R^d)$, smooth away from the origin and satisfying the
Marcinkiewicz-H\"{o}rmander-Mihlin condition
\footnote{$A\lesssim B$ simply means that there exists a universal constant 
$C>1$ so that $A\leq CB$. We will also sometime use the notation
$A\sim B$ to denote the statement that $A\lesssim B$ and $B\lesssim A$}

\begin{equation}\label{mhm}
|\partial^{\alpha} m(\xi) |\lesssim
\frac{1}{|\xi|^{|\alpha|}}
\end{equation}
for every $\xi\in \R^d\setminus\{0\}$ and sufficiently many multi-indices
$\alpha$. We say that such a symbol $m$ is {\it trivial} if and only if
$m(\xi)=1$ for every $\xi\in\R^d$.

If $n\geq 1$ is a fixed integer, we also denote by $\cal{M}$$_{flag}(\R^n)$
the set of all symbols $m$ given by arbitrary products of the form

\begin{equation}\label{mflag}
m(\xi):= \prod_{S\subseteq \{1,...,n\}}
m_S(\xi_S)
\end{equation}
where $m_S\in \cal{M}$$(\R^{card(S)})$, the vector $\xi_S\in\R^{card(S)}$ 
is defined by $\xi_S:= (\xi_i)_{i\in S}$, while $\xi\in\R^n$ is the vector 
$\xi:= (\xi_i)_{i=1}^n$.

Every symbol $m\in$$\cal{M}$$_{flag}(\R^n)$ defines an $n$-linear operator
$T_m$ by the formula

\begin{equation}\label{tmdef}
T_m(f_1, ..., f_n)(x) := \int_{\R^n}m(\xi)
\widehat{f_1}(\xi_1) ... 
\widehat{f_n}(\xi_n) e^{2\pi i x (\xi_1 + ... + \xi_n)} d\xi
\end{equation}
where $f_1, ..., f_n$ are Schwartz functions on the real line $\R$.

In the particular case when all the factors $(m_S)_{S\subseteq \{1,...,n\}}$
in (\ref{mflag}) are {\it trivial} the expression $T_m(f_1,...,f_n)(x)$ becomes
the product of our functions $f_1(x)\cdot...\cdot f_n(x)$ and as a consequence,
H\"{o}lder inequalities imply the fact that $T_m$ maps 
$L^{p_1}\times...\times L^{p_n}\rightarrow L^p$ boundedly as long as
$1<p_1,...,p_n <\infty$, $1/p_1+...+1/p_n = 1/p$ and $0<p<\infty$.
Similar estimates hold in the situation when all the factors 
$(m_S)_{S\subseteq \{1,...,n\}}$ in (\ref{mflag}) are {\it trivial} except for
the one corresponding to the set $\{1,...,n\}$. This deep and important fact
is a classical result in harmonic analysis known as the Coifman-Meyer theorem
\cite{meyerc}, \cite{gt}, \cite{ks}. 
Clearly, the same conclusion is also true if we assume that the only
{\it non-trivial} symbols are those corresponding to mutually disjoint subsets
of $\{1,...,n\}$, because this case can be factored out as a combination of
the previous two.

It is therefore natural to ask the following question.

\begin{question}\label{Q}
Is it true that $T_m$ maps 
$L^{p_1}\times...\times L^{p_n}\rightarrow L^p$ boundedly as long as
$1<p_1,...,p_n <\infty$, $1/p_1+...+1/p_n = 1/p$ and $0<p<\infty$ for any
$m\in$$\cal{M}$$_{flag}(\R^n)$ ?
\end{question}
The main goal of the present paper, is to give an affirmative answer to the
above question, in the simplest case which goes beyond the Coifman-Meyer 
theorem. We will consider the case of a tri-linear operator whose {\it non-trivial}
factors in (\ref{mflag}) are those corresponding to the subsets
$\{1,2\}$ and $\{2,3\}$.

More specifically, let $a, b\in$ $\cal{M}$$(\R^2)$ and denote by $T_{ab}$
the operator given by

\begin{equation}\label{tabdef}
T_{ab}(f_1, f_2, f_3)(x):= \int_{\R^3}
a(\xi_1, \xi_2) b(\xi_2, \xi_3)
\widehat{f_1}(\xi_1)
\widehat{f_2}(\xi_2)
\widehat{f_3}(\xi_3)
e^{2\pi i x(\xi_1 + \xi_2 + \xi_3)}
d\xi_1 d\xi_2 d\xi_3.
\end{equation}
Our main theorem is the following.

\begin{theorem}\label{main}
The operator $T_{ab}$ previously defined maps 
$L^{p_1}\times L^{p_2}\times L^{p_3}\rightarrow L^p$ as long as
$1<p_1, p_2, p_3 <\infty$, $1/p_1 + 1/p_2 + 1/p_3 = 1/p$ and
$0<p<\infty$.
\end{theorem}
Moreover, we will show that in this particular case there are also some
$L^{\infty}$-estimates available (in general, one cannot hope for any of them,
as one can easily see by taking all the factors in (\ref{mflag}) to be 
{\it trivial},
except for the ones corresponding to subsets which have cardinality $1$ ).
We believe however that the answer to our Question \ref{Q} is affirmative in general, and that
the $L^p$-estimates described above are satisfied by the operators
$T_m$ in (\ref{tmdef}) for all the symbols 
$m\in$$\cal{M}$$_{flag}(\R^n)$. We intend to address this general situation in 
a separate, future paper.

To motivate the introduction of these 
{\it paraproducts with flag singularities}, we should mention that some 
particular examples appeared implicitly in connection with the so-called
{\it bi-est} and {\it multi-est} operators studied in 
\cite{mtt:walshbiest}, \cite{mtt:fourierbiest}, \cite{mtt:multiest}.

The {\it bi-est} is the tri-linear operator $T_{bi-est}$ defined by the following
formula

\begin{equation}\label{biestdef}
T_{bi-est}(f_1, f_2, f_3)(x):=
\int_{\xi_1 < \xi_2 < \xi_3}
\widehat{f_1}(\xi_1)
\widehat{f_2}(\xi_2)
\widehat{f_3}(\xi_3)
e^{2\pi i x(\xi_1 + \xi_2 + \xi_3)}
d\xi_1 d\xi_2 d\xi_3
\end{equation}
and we know from \cite{mtt:walshbiest} and \cite{mtt:fourierbiest}
that it satisfies many $L^p$-estimates of the type described above.
Its symbol $\chi_{\xi_1 < \xi_2 < \xi_3}$ can be viewed as a product
of two bi-linear Hilbert transform type symbols, namely
$\chi_{\xi_1 < \xi_2}$ and $\chi_{\xi_2 < \xi_3}$ \cite{laceyt1}, 
\cite{laceyt2}. If one replaces them both with smoother symbols in the class
$\cal{M}$$(\R^2)$, then one obtains our tri-linear operator $T_{ab}$ in
(\ref{tabdef}).

As mentioned in \cite{mtt:bicarleson}, the interesting fact about
such operators as $T_m$ in (\ref{tmdef}), is that they have a very special
{\it multi-parameter} structure which seems to be new in harmonic analysis.
This structure is specific to the multi-linear analysis since only in this 
context one can construct operators given by multi-parameter symbols
which act on functions defined on the real line.

{\bf Acknowledgements:} The author has been partially supported by NSF and by
an Alfred P. Sloan Research Fellowship.

\section{Adjoint operators and interpolation}

The purpose of the present section is to recall the interpolation theory
from \cite{cct}, that will allow us to reduce our desired estimates in
Theorem \ref{main} to some restrictead weak type estimates, which are more
convenient to handle.

To each generic tri-linear operator $T$ we associate a four-linear form
$\Lambda$ defined by the following formula

\begin{equation}\label{form}
\Lambda (f_1, f_2, f_3, f_4):= \int_{\R} T(f_1, f_2, f_3)(x) f_4(x) dx.
\end{equation}

There are also three adjoint operators $T^{*j}$, $j=1,2,3$ attached to $T$,
defined by duality as follows

\begin{equation}\label{t*1}
\int_{\R} T^{*1}(f_2, f_3, f_4)(x) f_1(x) dx := \Lambda (f_1, f_2, f_3, f_4),
\end{equation}

\begin{equation}\label{t*2}
\int_{\R} T^{*2}(f_1, f_3, f_4)(x) f_2(x) dx := \Lambda (f_1, f_2, f_3, f_4),
\end{equation}

\begin{equation}\label{t*3}
\int_{\R} T^{*3}(f_1, f_2, f_4)(x) f_3(x) dx := \Lambda (f_1, f_2, f_3, f_4).
\end{equation}
For symmetry, we will also sometimes use the notation $T^{*4}:= T$. 

The following definition has been introduced in \cite{cct}.

\begin{definition}\label{weaktypedef}
Let $(p_1, p_2, p_3, p_4)$ be a $4$-tuple of real numbers so that
$1<p_1, p_2, p_3\leq \infty$, $1/p_1 + 1/p_2 + 1/p_3 = 1/p_4$ and
$0<p_4<\infty$. We say that the tri-linear operator $T$ is of restricted 
weak type $(p_1, p_2, p_3, p_4)$, if and only if for any $(E_i)_{i=1}^4$
measurable subsets of the real line $\R$ with $0< |E_i|<\infty$ for
$i=1,2,3,4$, there exists a subset $E'_4\subseteq E_4$ with
$|E'_4|\sim |E_4|$ so that

\begin{equation}\label{estimate}
|\int_{\R} T(f_1, f_2, f_3)(x) f_4(x) dx |
\lesssim 
|E_1|^{1/p_1}
|E_2|^{1/p_2}
|E_3|^{1/p_3}
|E_4|^{1/p'_4},
\end{equation}
for every $f_i\in X(E_i)$, $i=1,2,3$ and $f_4\in X(E'_4)$ where in general
$X(E)$ denotes the space of all measurable functions $f$ supported on $E$ with
$\|f\|_{\infty}\leq 1$ and $p'_4$ is the dual index of $p_4$ (note that
since $1/p_4 + 1/p'_4 = 1$, $p'_4$ can be negative if $0<p_4<1$ ).
\end{definition}

As in \cite{mtt:walshbiest}, \cite{mtt:fourierbiest} let us consider now the $3$ - dimensional
hyperspace $S$ defined by

$$S:= \{ (\alpha_1, \alpha_2, \alpha_3, \alpha_4)\in\R^4 : 
\alpha_1 + \alpha_2 + \alpha_3 + \alpha_4 = 1 \}.$$
Denote by $\bf{D}$ the open interior of the convex hull of the $7$ extremal points
$A_{11}$, 
$A_{12}$, 
$A_{21}$, 
$A_{22}$, 
$A_{31}$, 
$A_{32}$ and
$A_4$ in Figure $1$. They all belong to $S$ and have
the following coordinates: $A_{11} (-1, 1, 1, 0)$, 
$A_{12} (-1, 1, 0, 1)$,
$A_{21} (1, -1, 1, 0)$,
$A_{22} (0, 0, 0, 1)$,
$A_{31} (1, 1, -1, 0)$,
$A_{32} (0, 1, -1, 1)$ and
$A_4 (1, 1, 1, -2)$. 

\begin{figure}[htbp]\centering
\psfig{figure=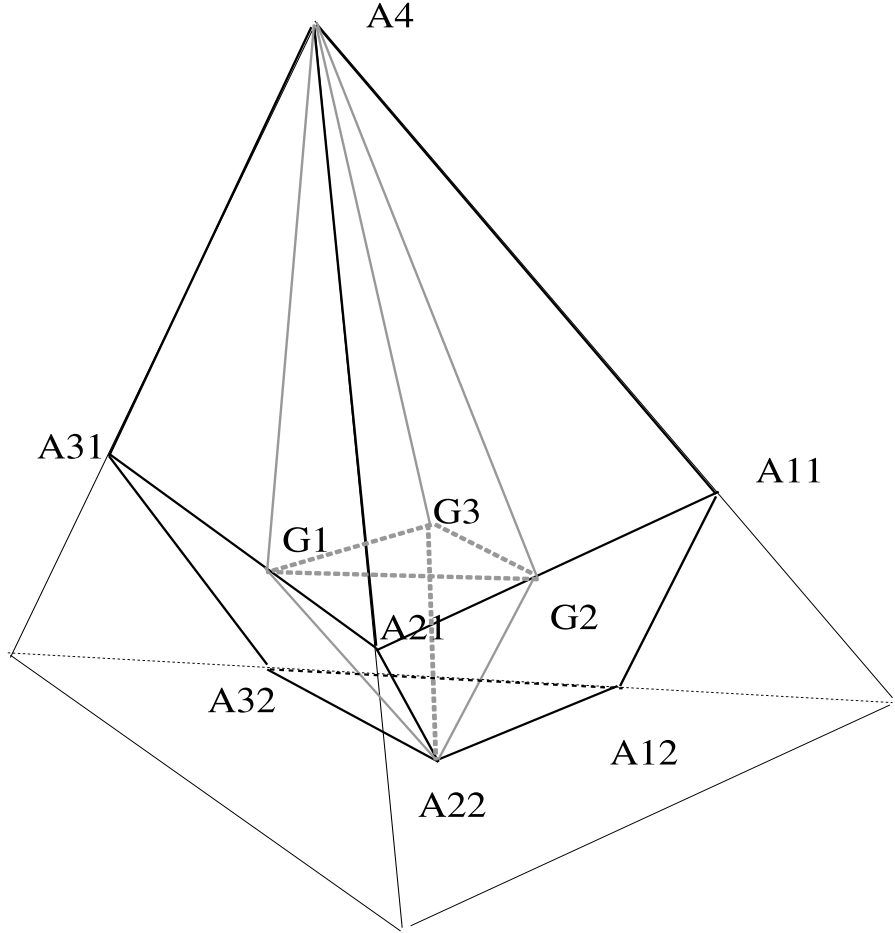, height=6in,width=5in}
\caption{Polytopes}
\label{fig}
\end{figure}

Denote also by $\widetilde{\bf{D}}$ the open interior of the convex hull of the $5$ extremal points
$A_{22}$, $G_1$, $G_2$, $G_3$ and $A_4$ where $G_1$, $G_2$ and $G_3$ have the coordinates
$(1, 0, 0, 0)$, $(0, 1, 0, 0)$ and $(0, 0, 1, 0)$ respectively. The following
theorem will be proved directly in the following sections.

\begin{theorem}\label{tabthm}
If $a, b\in\cal{M}$$(\R^2)$ are as before then, the following statements about the operator
$T_{ab}$ hold:

$(a)$ There exist points $(1/p_1, 1/p_2, 1/p_3, 1/p_4)\in\bf{D}$ arbitrarily close
to $A_4$ so that $T_{ab}$ is of restricted weak type $(p_1, p_2, p_3, p_4)$.

$(b)$ There exist points $(1/p^{ij}_1, 1/p^{ij}_2, 1/p^{ij}_3, 1/p^{ij}_4)\in\bf{D}$
arbitrarily close to $A_{ij}$ so that $T^{*i}_{ab}$ is of restricted weak type
$(p^{ij}_1, p^{ij}_2, p^{ij}_3, p^{ij}_4)$ for $i=1,2,3$ and $j=1,2$.
\end{theorem}
If we assume the above result, our main Theorem \ref{main} follows immediately from the 
interpolation theory developed in \cite{cct}. As a consequence of that theory, if
$(p_1, p_2, p_3, p)$ are so that $1< p_1, p_2, p_3\leq\infty$, 
$1/p_1 + 1/p_2 + 1/p_3 = 1/p$ and $0<p<\infty$ then $T_{ab}$ maps
$L^{p_1}\times L^{p_2}\times L^{p_3}\rightarrow L^p$ boundedly, as long as the point
$(1/p_1, 1/p_2, 1/p_3, 1/p)$ belongs to $\bf{D}$. And this is clearly true if
$(p_1, p_2, p_3, p)$ satisfies the hypothesis of Theorem \ref{main}.
In fact, in this case, the corresponding points $(1/p_1, 1/p_2, 1/p_3, 1/p)$
belong to $\widetilde{\bf{D}}$ which is a subset of $\bf{D}$.

Moreover, since we also observe that all the points of the form
$(0,\alpha, \beta,\gamma)$, $(\alpha, 0, \beta, \gamma)$, $(\alpha,\beta, 0,\gamma)$
with $\alpha, \beta >0$, $\alpha+\beta+\gamma=1$ and
$(0,\tilde{\alpha},0,\tilde{\beta})$ with $\tilde{\alpha}, \tilde{\beta} >0$ 
$\tilde{\alpha}+\tilde{\beta}=1$ belong to $\bf{D}$, we deduce that in addition
$T_{ab}$ maps $L^{\infty}\times L^p\times L^q\rightarrow L^r$,
$L^p\times L^{\infty}\times L^q\rightarrow L^r$, 
$L^p\times L^q\times L^{\infty}\rightarrow L^r$ and 
$L^{\infty}\times L^s\times L^{\infty}\rightarrow L^s$ boundedly, as long as
$1<p,q,s<\infty$, $0<r<\infty$ and $1/p+1/q=1/r$.

The only $L^{\infty}$ estimates that do not follow from such interpolation argumens are those of
the form $L^{\infty}\times L^{\infty}\times L^s\rightarrow L^s$ and
$L^s\times L^{\infty}\times L^{\infty}\rightarrow L^s$, because points of the form
$(1/s, 0,0,1/s')$ and $(0,0,1/s,1/s')$ only belong to the boundary of $\bf{D}$. But this is
not surprising since such estimates are false in general, as one can easily see by taking
$f_2\equiv1$ in (\ref{tabdef}).

In conclusion, to have a complete understanding of the boundedness properties of our
operator $T_{ab}$, it is enough to prove Theorem \ref{tabthm}.

\section{Discrete model operators}

In this section we introduce some discrete model operators and state a general theorem about 
them. Roughly speaking, this theorem says that they satisfy the desired restricted weak type
estimates in Theorem \ref{tabthm}. Later on, in Section $4$, we will prove that the analysis
of the operator $T_{ab}$ can in fact be reduced to the analysis of these discrete models.
We start with some notations.

An interval $I$ on the real line $\R$ is called dyadic if it is of the form
$I = [2^k n, 2^k (n+1)]$ for some $k,n\in\Z$. We denote by $\cal{D}$ the set of all such
dyadic intervals. If $J\in\cal{D}$ is fixed, we say that a smooth function
$\Phi_J$ is a bump adapted to $J$ if and only if the following inequalities hold

\begin{equation}\label{bump}
|\Phi_J^{(l)}(x)|\leq C_{l,\alpha}
\frac{1}{|J|^l}
\frac{1}{(1+\frac{\dist(x,J)}{|J|})^{\alpha}}
\end{equation}
for every integer $\alpha\in\N$ and sufficiently many derivatives $l\in\N$,
where $|J|$ is the length of $J$. If $\Phi_J$ is a bump adapted to $J$, we say
that $|J|^{-1/p}\Phi_J$ is an $L^p$- normalized bump adapted to $J$, for $1\leq p\leq \infty$.
We will also sometimes use the notation $\widetilde{\chi}_J$ for the approximate
cutoff function defined by

\begin{equation}\label{cutoff}
\widetilde{\chi}_J(x):= (1+\frac{\dist(x, J)}{|J|})^{-10}.
\end{equation}

\begin{definition}\label{non-lacunary}
A sequence of $L^2$- normalized bumps $(\Phi_I)_{I\in\cal{D}}$ adapted to dyadic intervals
$I\in\cal{D}$ is called a {\it non-lacunary} sequence if and only if for each
$I\in\cal{D}$ there exists an interval $\omega_I ( = \omega_{|I|})$  symmetric with respect to 
the origin
so that $\supp\widehat{\Phi_I}\subseteq \omega_I$ and $|\omega_I|\sim |I|^{-1}$.
\end{definition}

\begin{definition}\label{lacunary}
A sequence of $L^2$- normalized bumps $(\Phi_I)_{I\in\cal{D}}$ adapted to dyadic intervals
$I\in\cal{D}$ is called a {\it lacunary} sequence if and only if for each
$I\in\cal{D}$ there exists an interval $\omega_I ( = \omega_{|I|})$ 
so that $\supp\widehat{\Phi_I}\subseteq \omega_I$, 
$|\omega_I|\sim |I|^{-1}\sim\dist(0,\omega_I)$ and $0\notin 5\omega_I$.
\end{definition}
Let now consider $\I_1, \J_1\subseteq\D$ two finite families of dyadic intervals.
Let also $(\Phi^j_I)_{I\in\I_1}$ for $j=1,2,3$ be three sequences of $L^2$- normalized
bumps so that $(\Phi^2_I)_{I\in\I_1}$ is {\it non-lacunary} while
$(\Phi^j_I)_{I\in\I_1}$ for $j\neq 2$ are both {\it lacunary} in the sense of the above 
definitions. 

We also consider $(\Phi^j_J)_{J\in\J_1}$ for $j=1,2,3$ three sequences of $L^2$- normalized
bumps so that at least two of them are {\it lacunary}. Then, define the discrete model
operator $T_1$ by the formula

\begin{equation}\label{t1def}
T_1(f_1, f_2, f_3)(x):=
\sum_{I\in\I_1}
\frac{1}{|I|^{1/2}}
\langle f_1, \Phi^1_I \rangle
\langle B^1_I(f_2, f_3), \Phi^2_I \rangle
\Phi^3_I
\end{equation}
where

\begin{equation}\label{B1Idef}
B^1_I(f_2, f_3)(x):=
\sum_{J\in\J_1; |\omega^3_J|\leq |\omega^2_I|; \omega^3_J\cap \omega^2_I\neq\emptyset}
\frac{1}{|J|^{1/2}}
\langle f_2, \Phi^1_J \rangle
\langle f_3, \Phi^2_J \rangle
\Phi^3_J.
\end{equation}
If $k_0$ is a strictly positive integer, define also the operator $T_{1, k_0}$ by

\begin{equation}\label{t1k0def}
T_{1,k_0}(f_1, f_2, f_3)(x):=
\sum_{I\in\I_1}
\frac{1}{|I|^{1/2}}
\langle f_1, \Phi^1_I \rangle
\langle B^1_{I,k_0}(f_2, f_3), \Phi^2_I \rangle
\Phi^3_I
\end{equation}
where

\begin{equation}\label{B1Ik_0def}
B^1_{I,k_0}(f_2, f_3)(x):=
\sum_{J\in\J_1; 2^{k_0}|\omega^3_J|\sim |\omega^2_I|; \omega^3_J\cap \omega^2_I\neq\emptyset}
\frac{1}{|J|^{1/2}}
\langle f_2, \Phi^1_J \rangle
\langle f_3, \Phi^2_J \rangle
\Phi^3_J.
\end{equation}
Similarly, let us now consider two other finite families of dyadic intervals 
$\I_2, \J_2\subseteq \D$. As before, we also consider sequences
$(\Phi^j_I)_{I\in\I_2}$, $(\Phi^j_J)_{J\in\J_2}$ for $j=1,2,3$ of $L^2$- normalized bumps,
where this time we assume that $(\Phi^1_I)_{I\in\I_2}$ is {\it non-lacunary} while
$(\Phi^j_I)_{I\in\I_2}$ are both {\it lacunary} for $j\neq 1$ and at least two of the sequences
$(\Phi^j_J)_{J\in\J_2}$ are {\it lacunary}. Using them, we define the operator $T_2$ by
the formula

\begin{equation}\label{t2def}
T_2(f_1, f_2, f_3)(x):=
\sum_{I\in\I_2}
\frac{1}{|I|^{1/2}}
\langle B^2_I(f_1, f_2), \Phi^1_I \rangle
\langle f_3, \Phi^2_I \rangle
\Phi^3_I
\end{equation}
where

\begin{equation}\label{B2Idef}
B^2_I(f_1, f_2)(x):=
\sum_{J\in\J_2; |\omega^3_J|\leq |\omega^1_I|; \omega^3_J\cap \omega^1_I\neq\emptyset}
\frac{1}{|J|^{1/2}}
\langle f_1, \Phi^1_J \rangle
\langle f_2, \Phi^2_J \rangle
\Phi^3_J.
\end{equation}
And finally, as before, for any strictly positive integer $k_0$ define also the operator
$T_{2,k_0}$ by

\begin{equation}\label{t2k0def}
T_{2,k_0}(f_1, f_2, f_3)(x):=
\sum_{I\in\I_2}
\frac{1}{|I|^{1/2}}
\langle B^2_{I,k_0}(f_1, f_2), \Phi^1_I \rangle
\langle f_3, \Phi^2_I \rangle
\Phi^3_I
\end{equation}
where

\begin{equation}\label{B2Ik0def}
B^2_I(f_1, f_2)(x):=
\sum_{J\in\J_2; 2^{k_0}|\omega^3_J|\sim|\omega^1_I|; \omega^3_J\cap \omega^1_I\neq\emptyset}
\frac{1}{|J|^{1/2}}
\langle f_1, \Phi^1_J \rangle
\langle f_2, \Phi^2_J \rangle
\Phi^3_J.
\end{equation}
The following theorem about these operators will be proved carefully in the forthcoming sections.

\begin{theorem}\label{modelthm}
Our previous Theorem \ref{tabthm} holds also for all the operators $T_1$, $T_2$, $T_{1,k_0}$,
$T_{2, k_0}$ with bounds which are independent on $k_0$ and the cardinalities of the sets
$\I_1$, $\I_2$, $\J_1$, $\J_2$. Moreover, the subsets
$(E'_j)_{j=1}^4$ which appear implicitly due to Definition \ref{weaktypedef}, can be chosen
independently on the $L^2$- normalized families considered above.
\end{theorem}

\section{Reduction to the model operators}

As we promised, the aim of the present section is to show that the analysis of our operator
$T_{ab}$ can be indeed reduced to the analysis of the model operators defined in the previous
section. To achieve this, we will decompose the multipliers $a(\xi_1,\xi_2)$ and
$b(\xi_2, \xi_3)$ separetely and after that we will study their interactions.

Fix $M>0$ a big integer. For $j=1,2,...,M$ consider Schwartz functions
$\Psi_j$ so that $\supp\widehat{\Psi_j}\subseteq \frac{10}{9}[j-1, j]$,
$\Psi_j=1$ on $[j-1,j]$ and for $j=-M,-M+1,...,-1$ consider Schwartz functions
$\Psi_j$ so that $\supp\widehat{\Psi_j}\subseteq \frac{10}{9}[j, j+1]$
and $\Psi_j=1$ on $[j,j+1]$.
\footnote{If $I$ is an interval, we denote by $cI(c>0)$ the interval with the 
same center as $I$ and whose length is $c$ times the length of $I$}

If $\lambda$ is a positive real number and $\Psi$ is a Schwartz function, we 
denote by

$$D^p_{\lambda}\Psi(x):= \lambda^{-1/p}\Psi(\lambda^{-1}x)$$
the dilation operator which preserves the $L^p$ norm of $\Psi$, for
$1\leq p\leq\infty$.

Define the new symbol $\widetilde{a}(\xi_1, \xi_2)$ by the formula

\begin{equation}\label{atilde}
\widetilde{a}(\xi_1, \xi_2):=
\sum_{\max(|j'_1|, |j'_2|) =M}
\int_{\R}
D^{\infty}_{2^{\lambda'}}
\widehat{\Psi_{j'_1}}(\xi_1)
D^{\infty}_{2^{\lambda'}}
\widehat{\Psi_{j'_2}}(\xi_2)
d\lambda'.
\end{equation}
Clearly, by construction, $\widetilde{a}$ belongs to the class 
$\M(\R^2)$. Also, things can be arranged so that
$|\widetilde{a}(\xi_1, \xi_2)|\geq c_0>0$ for every $(\xi_1, \xi_2)\in\R^2$,
where $c_0$ is a universal constant. Roughly speaking, this 
$\widetilde{a}$ should be understood as being essentially a decomposition
of unity in frequency space into a series of smooth functions, supported on
rectangular annuli. Then, we write $a(\xi_1, \xi_2)$ as

$$a(\xi_1, \xi_2) = \frac{a(\xi_1, \xi_2)}{\widetilde{a}(\xi_1,\xi_2)}
\cdot\widetilde{a}(\xi_1,\xi_2):= \widetilde{\widetilde{a}}(\xi_1,\xi_2)
\cdot\widetilde{a}(\xi_1,\xi_2) = $$

\begin{equation}\label{1}
\sum_{\max(|j'_1|, |j'_2|) = M}
\int_{\R}
(\widetilde{\widetilde{a}}(\xi_1,\xi_2)
D^{\infty}_{2^{\lambda'}}
\widehat{\Psi_{j'_1}}(\xi_1)
D^{\infty}_{2^{\lambda'}}
\widehat{\Psi_{j'_2}}(\xi_2))
d\lambda'
\end{equation}
and observe that $\widetilde{\widetilde{a}}(\xi_1,\xi_2)$ has the same
properties as $a(\xi_1, \xi_2)$.

Fix now $j'_1, j'_2$ with $\max(|j'_1|, |j'_2|) = M$ and $\lambda'\in\R$. By taking
advantage of the fact that $\widetilde{\widetilde{a}}\in\M(\R^2)$, one can write it on the 
support of 
$D^{\infty}_{2^{\lambda'}}\widehat{\Psi_{j'_1}}\otimes
D^{\infty}_{2^{\lambda'}}\widehat{\Psi_{j'_2}}$
as a double Fourier series and this allows us to decompose the inner term in
(\ref{1}) as

$$\sum_{n'_1,n'_2\in\Z}
C(a)^{j'_1, j'_2}_{\lambda',n'_1, n'_2}
\left(D^{\infty}_{2^{\lambda'}}\widehat{\Psi_{j'_1}}(\xi_1)
e^{2\pi i n'_1\frac{9}{10}2^{-\lambda'}\xi_1}\right)
\left(D^{\infty}_{2^{\lambda'}}\widehat{\Psi_{j'_2}}(\xi_2)
e^{2\pi i n'_2\frac{9}{10}2^{-\lambda'}\xi_2}\right):=$$

$$\sum_{n'_1,n'_2\in\Z}
C(a)^{j'_1, j'_2}_{\lambda',n'_1, n'_2}
\left(D^{\infty}_{2^{\lambda'}}\widehat{\Psi^{n'_1}_{j'_1}}(\xi_1)\right)
\left(D^{\infty}_{2^{\lambda'}}\widehat{\Psi^{n'_2}_{j'_2}}(\xi_2)\right)$$
where we denoted by $\Psi^{n'_1}_{j'_1}$ and $\Psi^{n'_2}_{j'_2}$ the functions defined by

$$\widehat{\Psi^{n'_1}_{j'_1}}(\xi_1):= \widehat{\Psi_{j'_1}}(\xi_1)
e^{2\pi i n'_1\frac{9}{10}\xi_1}$$
and

$$\widehat{\Psi^{n'_2}_{j'_2}}(\xi_2):= \widehat{\Psi_{j'_2}}(\xi_2)
e^{2\pi i n'_2\frac{9}{10}\xi_2}$$
and the corresponding constants $C(a)^{j'_1, j'_2}_{\lambda',n'_1, n'_2}$ satisfy the 
inequalities

\begin{equation}\label{2}
|C(a)^{j'_1, j'_2}_{\lambda',n'_1, n'_2}|\lesssim
\frac{1}{(1+|n'_1|)^{1000}}
\frac{1}{(1+|n'_2|)^{1000}},
\end{equation}
for every $n'_1,n'_2\in\Z$, uniformly in $\lambda'\in\R$.

In particular, the symbol $a(\xi_1,\xi_2)$ can be written as

$$a(\xi_1,\xi_2) =$$

\begin{equation}\label{3}
\sum_{\max(|j'_1|, |j'_2|) = M}
\int_0^1
\sum_{n'_1,n'_2\in\Z}
\sum_{k'\in\Z}
C(a)^{j'_1, j'_2}_{k'+\k',n'_1, n'_2}
\left(D^{\infty}_{2^{k'+\k'}}\widehat{\Psi^{n'_1}_{j'_1}}(\xi_1)\right)
\left(D^{\infty}_{2^{k'+\k'}}\widehat{\Psi^{n'_2}_{j'_2}}(\xi_2)\right)
d\k'.
\end{equation}
Similarly, the symbol $b(\xi_2, \xi_3)$ can also be decomposed as

$$b(\xi_2,\xi_3)=$$

\begin{equation}\label{4}
\sum_{\max(|j''_1|, |j''_2|) = M}
\int_0^1
\sum_{n''_1,n''_2\in\Z}
\sum_{k''\in\Z}
C(b)^{j''_1, j''_2}_{k''+\k'',n''_1, n''_2}
\left(D^{\infty}_{2^{k''+\k''}}\widehat{\Psi^{n''_1}_{j''_1}}(\xi_2)\right)
\left(D^{\infty}_{2^{k''+\k''}}\widehat{\Psi^{n''_2}_{j''_2}}(\xi_3)\right)
d\k''
\end{equation}
where as before, the constants $C(b)^{j''_1, j''_2}_{k''+\k'',n''_1, n''_2}$
satisfy the inequalities

\begin{equation}\label{5}
|C(b)^{j''_1, j''_2}_{k''+\k'',n''_1, n''_2}|\lesssim
\frac{1}{(1+|n''_1|)^{1000}}
\frac{1}{(1+|n''_2|)^{1000}},
\end{equation}
for every $n''_1,n''_2\in\Z$, uniformly in $k''$ and $\k''$.

As a consequence, their product $a(\xi_1,\xi_2)\cdot b(\xi_2,\xi_3)$
becomes

$$a(\xi_1,\xi_2)\cdot b(\xi_2,\xi_3)= $$

$$
\sum_{\max(|j'_1|, |j'_2|) = M}
\sum_{\max(|j''_1|, |j''_2|) = M}
\sum_{n'_1,n'_2\in\Z}
\sum_{n''_1,n''_2\in\Z}\cdot
$$

$$\int_0^1 \int_0^1
\sum_{k',k''\in\Z}
C(a)^{j'_1, j'_2}_{k'+\k',n'_1, n'_2}\cdot
C(b)^{j''_1, j''_2}_{k''+\k'',n''_1, n''_2}\cdot
$$

\begin{equation}\label{6}
\left[
\left(D^{\infty}_{2^{k'+\k'}}\widehat{\Psi^{n'_1}_{j'_1}}(\xi_1)\right)
\left(D^{\infty}_{2^{k'+\k'}}\widehat{\Psi^{n'_2}_{j'_2}}(\xi_2)\right)\right]
\left[
\left(D^{\infty}_{2^{k''+\k''}}\widehat{\Psi^{n''_1}_{j''_1}}(\xi_2)\right)
\left(D^{\infty}_{2^{k''+\k''}}\widehat{\Psi^{n''_2}_{j''_2}}(\xi_3)\right)
\right]d\k' d\k''.
\end{equation}
Clearly, one has to have

\begin{equation}\label{7}
\supp (D^{\infty}_{2^{k'+\k'}}\widehat{\Psi^{n'_2}_{j'_2}})\cap
\supp (D^{\infty}_{2^{k''+\k''}}\widehat{\Psi^{n''_1}_{j''_1}})\neq\emptyset
\end{equation}
otherwise, the expression in (\ref{6}) vanishes.

Let now $\#$ be a positive integer, much bigger than $\log M$. 
If $k'$ and $k''$ are two integers
as in the sum above then, there are three possibilities:
either $k'\geq k''+\#$ or $k''\geq k'+\#$ or $|k'-k''|\leq \#$. As a consequence, the multiplier
$a(\xi_1,\xi_2)\cdot b(\xi_2,\xi_3)$ can be decomposed accordingly as

$$a(\xi_1,\xi_2)\cdot b(\xi_2,\xi_3) =
m_1(\xi_1,\xi_2,\xi_3) + m_2(\xi_1,\xi_2,\xi_3) + m_3(\xi_1,\xi_2,\xi_3).$$
Since it is not difficult to see that $m_3\in\M(\R^2)$, the desired estimates
for the tri-linear operator $T_{m_3}$ follow from the classical Coifman-Meyer theorem quoted 
before. It is therefore enough to concentrate our attention on the remaining operators
$T_{m_1}$ and $T_{m_2}$. Since their definitions are symmetric, we will only study the case
of $T_{m_1}$ where the summation in (\ref{6}) runs over those
$k', k''\in\Z$ having the property that $k'\geq k''+\#$. We then observe that
since $\#$ is big in comparison to $\log M$, we have to have $j'_2=-1$ or
$j'_2=1$ in order for (\ref{7}) to hold. In particular, this simplies that
the intervals $\supp (D^{\infty}_{2^{k'+\k'}}\widehat{\Psi^{n'_2}_{j'_2}})_{k'\in\Z}$
are all intersecting each other.

At this moment, let us also remind ourselves that in order to prove restricted weak
type estimates for $T_{m_1}$, we would need to understand expressions of the form

$$\left|
\int_{\R}T_{m_1}(f_1, f_2, f_3)(x) f_4(x) dx\right| =
$$

\begin{equation}\label{77}
\left|
\int_{\xi_1+\xi_2+\xi_3+\xi_4=0}
m_1(\xi_1, \xi_2, \xi_3)
\widehat{f_1}(\xi_1)
\widehat{f_2}(\xi_2)
\widehat{f_3}(\xi_3)
\widehat{f_4}(\xi_4)d\xi\right|
\end{equation}
and as a consequence, from now on, we will think of our tri-dimensional vectors
$(\xi_1,\xi_2, \xi_3)\in\R^3$ as being part of $4$- dimensional ones
$(\xi_1, \xi_2, \xi_3, \xi_4)\in\R^4$ for which $\xi_1+\xi_2+\xi_3+\xi_4=0$.

Fix now the parameters $j'_1$, $j'_2$, $j''_1$, $j''_2$, $n'_1$, $n'_2$, $n''_1$,
$n''_2$, $k'$, $k''$, $\k'$, $\k''$ so that 
$k'\geq k''+\#$ and look at the corresponding inner term in (\ref{6}). It can be rewritten as

\begin{equation}\label{8}
\widehat{\Psi^{n'_1}_{k',\k',j'_1}}(\xi_1)
\widehat{\Psi^{n'_2}_{k',\k',j'_2}}(\xi_2)
\widehat{\Psi^{n''_1}_{k'',\k'',j''_1}}(\xi_2)
\widehat{\Psi^{n''_2}_{k'',\k'',j''_2}}(\xi_3)
\end{equation}
where

$$\Psi^{n'_1}_{k',\k',j'_1}:= D^1_{2^{-k'-\k'}}\Psi^{n'_1}_{j'_1},$$

$$\Psi^{n'_2}_{k',\k',j'_2}:= D^1_{2^{-k'-\k'}}\Psi^{n'_2}_{j'_2},$$

$$\Psi^{n''_1}_{k'',\k'',j''_1}:= D^1_{2^{-k''-\k''}}\Psi^{n''_1}_{j''_1}$$
and

$$\Psi^{n''_2}_{k'',\k'',j''_1}:= D^1_{2^{-k''-\k''}}\Psi^{n''_2}_{j''_2}.$$
Consider now Schwartz functions $\Psi_{k',\k',j'_1,j'_2}$ and
$\Psi_{k'',\k'',j''_1,j''_2}$ so that 
$\widehat{\Psi_{k',\k',j'_1,j'_2}}$ is identically equal to $1$ on the interval
$-2(\supp (\widehat{\Psi^{n'_1}_{k',\k',j'_1}}) + 
\supp (\widehat{\Psi^{n'_2}_{k',\k',j'_2}}))$ and is supported on a $\frac{10}{9}$
enlargement of it, while $\widehat{\Psi_{k'',\k'',j''_1,j''_2}}$ is identically
equal to $1$ on the interval
$(\supp (\widehat{\Psi^{n''_1}_{k'',\k'',j''_1}}) + 
\supp (\widehat{\Psi^{n''_2}_{k'',\k'',j''_2}}))$ and is also supported on a
$\frac{10}{9}$ enlargement of it.
Since $\xi_1+\xi_2+\xi_3+\xi_4=0$, one can clearly insert these two new functions into the
previous expression (\ref{8}), which now becomes

\begin{equation}\label{9}
\left[
\widehat{\Psi^{n'_1}_{k',\k',j'_1}}(\xi_1)
\widehat{\Psi^{n'_2}_{k',\k',j'_2}}(\xi_2)
\widehat{\Psi_{k',\k',j'_1,j'_2}}(\xi_4)\right]\cdot
\end{equation}

$$
\left[
\widehat{\Psi^{n''_1}_{k'',\k'',j''_1}}(\xi_2)
\widehat{\Psi^{n''_2}_{k'',\k'',j''_2}}(\xi_3)
\widehat{\Psi_{k'',\k'',j''_1,j''_2}}(\xi_2+\xi_3)\right].
$$
The following elementary lemmas will play an important role in our further 
decomposition (see also \cite{mtt:fourierbiest}).

\begin{lemma}\label{calc1}
Let $\eta_1, \eta_2,\eta_3,\eta_4,\eta_{14},\eta_{23}$ be Schwartz functions.
Then,

$$\int_{\xi_1+\xi_2+\xi_3+\xi_4=0}
\widehat{\eta_1}(\xi_1)
\widehat{\eta_2}(\xi_2)
\widehat{\eta_3}(\xi_3)
\widehat{\eta_4}(\xi_4)
\widehat{\eta_{14}}(\xi_1+\xi_4)
\widehat{\eta_{23}}(\xi_2+\xi_3)
\widehat{f_1}(\xi_1)
\widehat{f_2}(\xi_2)
\widehat{f_3}(\xi_3)
\widehat{f_4}(\xi_4)
d\xi =$$

$$\int_{\R}
\left[(f_1*\eta_1)(f_4*\eta_4)\right]*\eta_{14}\cdot
\left[(f_2*\eta_2)(f_3*\eta_3)\right]*\eta_{23}
dx.$$
\end{lemma}

\begin{proof}
We write the left hand side of the identity as

$$\int_{\xi_1+\xi_2+\xi_3+\xi_4=0}
\widehat{f_1*\eta_1}(\xi_1)
\widehat{f_4*\eta_4}(\xi_4)
\widehat{\eta_{14}}(\xi_1+\xi_4)
\widehat{f_2*\eta_2}(\xi_2)
\widehat{f_3*\eta_3}(\xi_3)
\widehat{\eta_{23}}(\xi_2+\xi_3)
d\xi =
$$

$$
\int_{\R}
\left[
\int_{\xi_1+\xi_4=\lambda}
\widehat{f_1*\eta_1}(\xi_1)
\widehat{f_4*\eta_4}(\xi_4)d\xi_1 d\xi_4\right]
\widehat{\eta_{14}}(\lambda)\cdot
$$

$$
\left[
\int_{\xi_2+\xi_3=-\lambda}
\widehat{f_2*\eta_2}(\xi_2)
\widehat{f_3*\eta_3}(\xi_3)d\xi_2 d\xi_3\right]
\widehat{\eta_{23}}(-\lambda)
d\lambda =$$

$$\int_{\R}
\left[
\widehat{(f_1*\eta_1)(f_4*\eta_4)}(\lambda)\widehat{\eta_{14}}(\lambda)\right]\cdot
$$

$$\left[
\widehat{(f_2*\eta_2)(f_3*\eta_3)}(-\lambda)\widehat{\eta_{23}}(-\lambda)\right]
d\lambda =
$$

$$\int_{\R}
\widehat{\left[(f_1*\eta_1)(f_4*\eta_4) \right]*\eta_{14}}(\lambda)\cdot
\widehat{\left[(f_2*\eta_2)(f_3*\eta_3) \right]*\eta_{23}}(-\lambda)
d\lambda
$$
and this, by Plancherel, is equal to the right hand side of the identity.

\end{proof}

\begin{lemma}\label{calc2}
Let $k\in\Z$ be a fixed integer, $F_1$, $F_2$, $F_3$ three functions
in $L^1\cap L^{\infty}(\R)$ and $\Phi_1$, $\Phi_2$, $\Phi_3$ three
$L^1$ normalized bumps adapted to the interval $[0,2^k]$.
Then,

\begin{equation}\label{10}
\int_{\R}
(F_1*\Phi_1)(x)
(F_2*\Phi_2)(x)
(F_3*\Phi_3)(x)
dx =
\end{equation}

$$\int_0^1
\sum_{I\in\D; |I|=2^k}
\frac{1}{|I|^{1/2}}
\langle F_1, \Phi_{I, t, 1}\rangle
\langle F_2, \Phi_{I, t, 2}\rangle
\langle F_3, \Phi_{I, t, 3}\rangle
dt$$
where $\Phi_{I, t, j}(y):= |I|^{1/2}\overline{F_j(x_I + t|I| -y)}$
for $j=1,2,3$ and $x_I$ is the left hand side of the dyadic interval $I$.
\end{lemma}

\begin{proof} For every $j=1,2,3$ write

$$F_j*\Phi_j(x) = 
\int_{\R}F_j(y)\Phi_j(x-y) dy
= 2^{-k/2}
\langle F_j, \Phi_{x,j}\rangle$$
where $\Phi_{x, j}(y):= 2^{k/2}\overline{\Phi_j(x-y)}$.

As a consequence, the left hand side of (\ref{10}) becomes

$$2^{-3k/2}
\int_{\R}
\langle F_1, \Phi_{x, 1}\rangle
\langle F_2, \Phi_{x, 2}\rangle
\langle F_3, \Phi_{x, 3}\rangle
dx =$$

$$2^{-3k/2}
\sum_{I\in\D; |I|=2^k}
\int_I
\langle F_1, \Phi_{x, 1}\rangle
\langle F_2, \Phi_{x, 2}\rangle
\langle F_3, \Phi_{x, 3}\rangle
dx =$$

$$2^{-3k/2}
\sum_{I\in\D; |I|=2^k}
\int_0^{2^k}
\langle F_1, \Phi_{x_I+z, 1}\rangle
\langle F_2, \Phi_{x_I+z, 2}\rangle
\langle F_3, \Phi_{x_I+z, 3}\rangle
dz.
$$
If we now change the variables by writing $z=t|I|$, then this expression becomes
precisely the right hand side of (\ref{10}).

\end{proof}

As a consequence, we have the following corollary.

\begin{corollary}\label{calc3}
Let $k', k''\in\Z$ be as before, $\Psi_1$, $\Psi_4$, $\Psi_{14}$ be three
$L^1$ normalized bumps adapted to the interval $[0, 2^{-k'}]$ and $\Psi_2$,
$\Psi_3$, $\Psi_{23}$ be three bumps adapted to the interval $[0, 2^{k''}]$.Then,

\begin{equation}\label{11}
\int_{\xi_1+\xi_2+\xi_3+\xi_4=0}
\widehat{\Psi_1}(\xi_1)
\widehat{\Psi_2}(\xi_2)
\widehat{\Psi_3}(\xi_3)
\widehat{\Psi_4}(\xi_4)
\widehat{\Psi_{14}}(\xi_1+\xi_4)
\widehat{\Psi_{23}}(\xi_2+\xi_3)
\widehat{f_1}(\xi_1)
\widehat{f_2}(\xi_2)
\widehat{f_3}(\xi_3)
\widehat{f_4}(\xi_4)
d\xi =
\end{equation}

$$\int_0^1
\sum_{I\in\D; |I|=2^{-k'}}
\langle f_1, \Psi_{I, t', 1}\rangle
\langle B_{k''}(f_2, f_3), \widetilde{\Psi}_{I, t', 14}\rangle
\langle f_4, \Psi_{I, t', 4}\rangle
d t'$$
where $B_{k''}(f_2, f_3)$ is given by

$$B_{k''}(f_2, f_3)(x) =
\int_0^1
\sum_{J\in\D; |J|=2^{-k''}}
\langle f_2, \Psi_{J, t'', 2}\rangle
\langle f_3, \Psi_{J, t'', 3}\rangle
\overline{\widetilde{\Psi}}_{J, t'', 23}(x)
d t''
$$
while $\widetilde{\Psi}_{I, t', 14}(y):= |I|^{1/2}\overline{\Psi_{14}(y- x_I -t'|I|)}$
and $\widetilde{\Psi}_{J, t'', 23}(y):= |J|^{1/2}\overline{\Psi_{23}(y- x_J- t''|J|)}$.
\end{corollary}

\begin{proof}
By using the first Lemma \ref{calc1}, the left hand side of (\ref{11}) is equal to

$$\int_{\R}
\left[(f_1*\Psi_1)(f_4*\Psi_4)\right]*\Psi_{14}(x)
\left[(f_2*\Psi_2)(f_3*\Psi_3)\right]*\Psi_{23}(x)
dx =$$

$$\int_{\R}
\left[(f_1*\Psi_1)(f_4*\Psi_4)\right]*\Psi_{14}(x) F_{23}(x) dx =$$

$$\int_{\R}
(f_1*\Psi_1)(x)(f_4*\Psi_4)(x)(F_{23}*\widetilde{\Psi}_{14})(x) dx,
$$
where $\widetilde{\Psi}_{14}$ is the reflection of $\Psi_{14}$ defined by
$\widetilde{\Psi}_{14}(y):= \Psi_{14}(-y)$ and $F_{23}$ is given by

$$F_{23}(x):= \left[(f_2*\Psi_2)(f_3*\Psi_3)\right]*\Psi_{23}(x).$$
By using the second Lemma \ref{calc2}, this can be further decomposed as

$$\int_0^1
\sum_{I\in\D; |I|=2^{-k'}}
\langle f_1, \Psi_{I, t', 1}\rangle
\langle F_{23}, \widetilde{\Psi}_{I, t', 14}\rangle
\langle f_4, \Psi_{I, t', 4}\rangle
d t'$$
On the other hand, since $\langle F_{23}, \widetilde{\Psi}_{I, t', 14}\rangle$ can also be 
written as

$$\int_{\R}
F_{23}(x)\overline{\widetilde{\Psi}}_{I, t', 14}(x) dx
= \int_{\R}\left[(f_2*\Psi_2)(f_3*\Psi_3)\right]*\Psi_{23}(x)
\overline{\widetilde{\Psi}}_{I, t', 14}(x) dx =$$

$$\int_{\R}
(f_2*\Psi_2)(x)(f_3*\Psi_3)(x)(\overline{\widetilde{\Psi}}_{I, t', 14}*\tilde{\Psi}_{23})(x)
dx,
$$
we can apply again Lemma \ref{calc1} and this will lead us to the desired expression.
\end{proof}

Clearly, modulo the two averages over parameters $t', t''\in [0,1]$, the
discretized expressions in Corollary \ref{calc3} are similar to the ones that appeared
in the definition of the model operators $T_1$ and $T_{1, k_0}$ in Section $3$ 
(one has to consider
the $4$- linear form associated to them to have a perfect similarity). Consequently,
we would like to apply this corollary to the expressions obtained after combining (\ref{9})
with (\ref{77}). We observe however that the formulas in (\ref{9}) are not precisely
of the required form (we would need to have instead of the factor
$\widehat{\Psi^{n_2'}_{k',\k',j'_2}}(\xi_2)$ a similar one but depending on $\xi_1+\xi_4$)
and so they need to be ``fixed''.

Before doing this, let us first make another reduction. Write the operator $T_{m_1}$ as

\begin{equation}\label{12}
T_{m_1}:=
\sum_{\max(|j'_1|, |j'_2|) = M}
\sum_{\max(|j''_1|, |j''_2|) = M}
\sum_{n'_1, n'_2\in\Z}
\sum_{n''_1, n''_2\in\Z}
T^{\vec{j'},\vec{j''}}_{\vec{n'},\vec{n''}},
\end{equation}
where $T^{\vec{j'},\vec{j''}}_{\vec{n'},\vec{n''}}$ are given by the correspoding
symbols in (\ref{6}) with the expressions in (\ref{6}) being replaced by their new formulas in
(\ref{9}) and where the summation over $k', k''\in\Z$ satisfies the constraint
$k'\geq k''+\#$.

We are going to prove explicitly that for each $\vec{j'}:= (j'_1, j'_2)$ and 
$\vec{j''}:= (j''_1, j''_2)$ the operator
$T^{\vec{j'},\vec{j''}}_{\vec{0},\vec{0}}:= T^{\vec{j'},\vec{j''}}$ satisfies the required
estimates. It will also be clear from our proof that the same arguments give

$$
\left\|
T^{\vec{j'},\vec{j''}}_{\vec{n'},\vec{n''}}
\right\|_{L^{p_1}\times L^{p_2}\times L^{p_3}\rightarrow L^p}
\lesssim
\frac{1}{(1+|\vec{n'}|)^{10}}
\frac{1}{(1+|\vec{n''}|)^{10}}
\left\|
T^{\vec{j'},\vec{j''}}
\right\|_{L^{p_1}\times L^{p_2}\times L^{p_3}\rightarrow L^p},
$$
and this would be enough to prove our desired estimates for $T_{m_1}$, due to the big decay in
(\ref{2}).

We now come back to the operator $T^{\vec{j'},\vec{j''}}$. Its symbol is given by an infinite
sum of expressions of the form (see (\ref{9}) and (\ref{12}))

\begin{equation}\label{13}
\left[
\widehat{\Psi_{k',\k',j'_1}}(\xi_1)
\widehat{\Psi_{k',\k',j'_2}}(\xi_2)
\widehat{\Psi_{k',\k',j'_1,j'_2}}(\xi_4)\right]\cdot
\end{equation}

$$
\left[
\widehat{\Psi_{k'',\k'',j''_1}}(\xi_2)
\widehat{\Psi_{k'',\k'',j''_2}}(\xi_3)
\widehat{\Psi_{k'',\k'',j''_1,j''_2}}(\xi_2+\xi_3)\right],
$$
where we suppressed the indices $n'_1,n'_2, n''_1, n''_2$, since they are all equal to zero now.

Fix then $\widetilde{M}\in [100, 200]$ an integer and write the function
$\widehat{\Psi_{k',\k',j'_2}}(\xi_2)$ as a Taylor series as follows

$$\widehat{\Psi_{k',\k',j'_2}}(\xi_2) = 
\sum_{l=0}^{\widetilde{M}-1}
(-\xi_3)^l
\frac{(\widehat{\Psi_{k',\k',j'_2}})^{(l)}(\xi_2+\xi_3)}{l!}
+ (-\xi_3)^{\widetilde{M}}
\frac{R^{\widetilde{M}}_{k',\k',j'_2}(\xi_2, \xi_3)}{\widetilde{M}!} =
$$

$$\sum_{l=0}^{\widetilde{M}-1}
\frac{(-\xi_3)^l}{l!}
(\widehat{\Psi_{k',\k',j'_2}})^{(l)}(-\xi_1 -\xi_4) +
\frac{(-\xi_3)^{\widetilde{M}}}{\widetilde{M}!}
R^{\widetilde{M}}_{k',\k',j'_2}(\xi_2, \xi_3),
$$
where $R^{\widetilde{M}}_{k',\k',j'_2}$ is the usual $\widetilde{M}$th rest in the Taylor
expansion.

Inserting this into (\ref{13}) we rewrite (\ref{13}) as

\begin{equation}\label{14}
\sum_{l=0}^{\widetilde{M}-1}
(\frac{2^{k''}}{2^{k'}})^l
\left[
\widehat{\Psi_{k',\k',j'_1}}(\xi_1)
\widehat{\Psi_{k',\k',j'_2,l}}(\xi_1+\xi_4)
\widehat{\Psi_{k',\k',j'_1,j'_2}}(\xi_4)\right]\cdot
\end{equation}

$$\left[
\widehat{\Psi_{k'',\k'',j''_1}}(\xi_2)
\widehat{\Psi_{k'',\k'',j''_2,l}}(\xi_3)
\widehat{\Psi_{k'',\k'',j''_1,j''_2}}(\xi_2+\xi_3)\right] +
$$

$$(\frac{2^{k''}}{2^{k'}})^{\widetilde{M}}
m_{\vec{k}, \vec{\k},\vec{j'},\vec{j''},\widetilde{M}}(\xi_1, \xi_2, \xi_3, \xi_4),
$$
where the functions $\widehat{\Psi_{k',\k',j'_2,l}}(\xi_1+\xi_4)$,
$\widehat{\Psi_{k'',\k'',j''_2,l}}(\xi_3)$ and
$m_{\vec{k}, \vec{\k},\vec{j'},\vec{j''},\widetilde{M}}(\xi_1, \xi_2, \xi_3, \xi_4)$ have
the obvious definitions, so that the two expressions in (\ref{13}) and (\ref{14})
to be consistent.

In particular, using (\ref{12}) and (\ref{14}) one can decompose 
$T^{\vec{j'},\vec{j''}}$ accordingly as

$$T^{\vec{j'},\vec{j''}} = T^{\vec{j'},\vec{j''}}_0
+ \sum_{l=1}^{\widetilde{M}-1}T^{\vec{j'},\vec{j''}}_l +
T^{\vec{j'},\vec{j''}}_{\widetilde{M}}.$$
Since we are in the case when $k'\geq k''+\#$, we can decompose 
$T^{\vec{j'},\vec{j''}}$ even further as

\begin{equation}\label{15}
T^{\vec{j'},\vec{j''}} = T^{\vec{j'},\vec{j''}}_0 +
\sum_{l=1}^{\widetilde{M}-1}\sum_{k_0=\#}^{\infty}
(2^{-k_0})^lT^{\vec{j'},\vec{j''}}_{l,k_0} +
\sum_{k_0=\#}^{\infty}
(2^{-k_0})^{\widetilde{M}}
T^{\vec{j'},\vec{j''}}_{\widetilde{M}, k_0}.
\end{equation}
For a fixed $k_0\geq \#$ we observe that the multiplier corresponding to the operator
$T^{\vec{j'},\vec{j''}}_{\widetilde{M}, k_0}$
which we denote by
$m_{\vec{j'},\vec{j''},\widetilde{M},k_0}(\xi_1, \xi_2, \xi_3)$ satisfies the
estimates

$$|\partial^{\alpha}
m_{\vec{j'},\vec{j''},\widetilde{M},k_0}(\vec{\xi})|
\lesssim (2^{k_0})^{|\alpha|}\frac{1}{|\vec{\xi}|^{|\alpha|}}
$$
for sufficiently many multi-indices $\alpha$ and as a consequence the classical
Coifman-Meyer theorem (see for instance its new proof in \cite{mptt:biparameter})
provides the required estimates for $T^{\vec{j'},\vec{j''}}_{\widetilde{M}, k_0}$
with a bound not bigger than $C2^{10 k_0}$, which is acceptable due
to the big decay in (\ref{15}). It is theorefore enough to understand the operators
$T^{\vec{j'},\vec{j''}}_0$ and $T^{\vec{j'},\vec{j''}}_{l,k_0}$ for $l=1,...,\widetilde{M}-1$
and $k_0\geq \#$. But their multipliers have the correct form now and to them
we can apply the discretization procedure provided by Corollary \ref{calc3}. And this will reduce
them to the model operators $T_1$ and $T_{1,k_0}$ defined in Section $3$.
\footnote{The {\it lacunarity} and {\it non-lacunarity} assumptions are also satisfied, 
as one can easily check}
Using now Theorem \ref{modelthm} and tacking advantage of the uniformity properties
described there, the estimates for $T^{\vec{j'},\vec{j''}}$ follow imediately.

In conclusion, it is indeed sufficient to prove our estimates for these model operators.

\section{$L^{1,\infty}$- sizes and $L^{1,\infty}$- energies}

We can now start the proof of Theorem \ref{modelthm}. It is of course enough to treat the 
operators $T_1$ and $T_{1,k_0}$ only, since the case of $T_2$ and $T_{2,k_0}$
is similar. We denote by $\Lambda_1$ and $\Lambda_{1,k_0}$ the $4$- linear forms associated
with the operators $T_1$ and $T_{1,k_0}$. As in \cite{mtt:fourierbiest}, 
since the $I$- spatial intervals are narrower
than their corresponding $J$- spatial intervals, it will be convenient to change the order
of summation in (\ref{t1def}) and rewrite the form $\Lambda_1$ as

\begin{equation}\label{form1}
\Lambda_1(f_1, f_2, f_3, f_4) =
\sum_{J\in\J_1}
\frac{1}{|J|^{1/2}}
a_J^{(1)}
a_J^{(2)}
a_J^{(3)}
\end{equation}
where

$$a_J^{(1)}:= \langle f_2, \Phi_J^1\rangle$$

$$a_J^{(2)}:= \langle f_3, \Phi_J^2\rangle$$
and

$$a_J^{(3)}:=
\langle
\sum_{I\in\I_1; \omega_J^3\cap \omega_I^2\neq\emptyset; |\omega_J^3|\leq |\omega_I^2|}
\frac{1}{|I|^{1/2}}
\langle f_1, \Phi_I^1\rangle
\langle f_4, \Phi_I^3\rangle
\Phi_I^2, \Phi_J^3\rangle.
$$
Similarly, we rewrite the form $\Lambda_{1, k_0}$ as

\begin{equation}\label{form1k0}
\Lambda_{1,k_0}(f_1, f_2, f_3, f_4) =
\sum_{J\in\J_1}
\frac{1}{|J|^{1/2}}
a_J^{(1)}
a_J^{(2)}
a_{J,k_0}^{(3)}
\end{equation}
where

$$a_{J,k_0}^{(3)}:=
\langle
\sum_{I\in\I_1; \omega_J^3\cap \omega_I^2\neq\emptyset; 2^{k_0}|\omega_J^3|\sim |\omega_I^2|}
\frac{1}{|I|^{1/2}}
\langle f_1, \Phi_I^1\rangle
\langle f_4, \Phi_I^3\rangle
\Phi_I^2, \Phi_J^3\rangle.
$$
We know from the definition of $T_1$ and $T_{1,k_0}$ in Section $3$ that the family
$(\Phi_I^2)_I$ may be {\it non-lacunary} while $(\Phi_I^i)_I$ for $i\neq 2$ are both 
{\it lacunary}. On the other hand we also know that there exists a unique $j=1,2,3$
which we fix from now on, so that the corresponding family
$(\Phi_J^j)_J$ is {\it non-lacunary} while $(\Phi_J^i)_J$ for $i\neq j$ are both
{\it lacunary}.

The standard way to estimate the forms $\Lambda_1$ and $\Lambda_{1,k_0}$ is to do so
by introducing some {\it sizes} and {\it energies} which in our case are going to be
more abstract variants of similar quantities considered in \cite{mptt:biparameter}.

The following definition contains those expressions which will be useful when estimating
the form $\Lambda_1$.

\begin{definition}\label{se1}
Let $\J$ be a finite family of dyadic intervals and $i=1,2,3$. For $i=j$ we define

$$\size^j_{i,\J}((a_J^{(i)})_{J}):=
\sup_{J\in\J}\frac{|a_J^{(i)}|}{|J|^{1/2}}
$$
and for $i\neq j$ we define

$$\size^j_{i,\J}((a_J^{(i)})_{J}):=
\sup_{J\in\J}
\frac{1}{|J|}
\left\|\left(
\sum_{J'\in\J; J'\subseteq J}
\frac{|a_{J'}^{(i)}|^2}{|J'|}\chi_{J'}(x)
\right)^{1/2}\right\|_{1,\infty}.
$$
Similarly, for $i= j$ we define

$$\energy^j_{i,\J}((a_J^{(i)})_{J}):=
\sup_{n\in\Z}\sup_{\bf{D}} 2^n (\sum_{J\in\bf{D}} |J| )
$$
where $\bf{D}$ ranges over those collections of disjoint dyadic intervals $J$ having the property
that

$$\frac{|a_J^{(i)}|}{|J|^{1/2}}\geq 2^n$$
and finally, for $i\neq j$ we define

$$\energy^j_{i,\J}((a_J^{(i)})_{J}):=
\sup_{n\in\Z}\sup_{\bf{D}} 2^n (\sum_{J\in\bf{D}} |J| )
$$
where this time $\bf{D}$ ranges over those collections of disjoint dyadic intervals $J$ having 
the property that

$$\frac{1}{|J|}
\left\|\left(
\sum_{J'\in\J; J'\subseteq J}
\frac{|a_{J'}^{(i)}|^2}{|J'|}\chi_{J'}(x)
\right)^{1/2}\right\|_{1,\infty}\geq 2^n.$$
\end{definition}

The next definition will be useful when estimating the form $\Lambda_{1,k_0}$.

\begin{definition}\label{se1}
Let $\J$ be a finite family of dyadic intervals and $k_0\geq \#$. 
For $j=3$ we define

$$\size^j_{3,k_0,\J}((a_{J,k_0}^{(3)})_{J}):=
\sup_{J\in\J}\frac{|a_{J,k_0}^{(3)}|}{|J|^{1/2}}
$$
and for $j\neq 3$ we define

$$\size^j_{3,k_0,\J}((a_{J,k_0}^{(3)})_{J}):=
\sup_{J\in\J}
\frac{1}{|J|}
\left\|\left(
\sum_{J'\in\J; J'\subseteq J}
\frac{|a_{J',k_0}^{(3)}|^2}{|J'|}\chi_{J'}(x)
\right)^{1/2}\right\|_{1,\infty}.
$$
Similarly, for $j=3$ we define

$$\energy^j_{3,k_0,\J}((a_{J,k_0}^{(3)})_{J}):=
\sup_{n\in\Z}\sup_{\bf{D}} 2^n (\sum_{J\in\bf{D}} |J| )
$$
where $\bf{D}$ ranges over those collections of disjoint dyadic intervals $J$ having the property
that

$$\frac{|a_{J,k_0}^{(3)}|}{|J|^{1/2}}\geq 2^n$$
and finally, for $j\neq 3$ we define

$$\energy^j_{3,k_0,\J}((a_{J,k_0}^{(3)})_{J}):=
\sup_{n\in\Z}\sup_{\bf{D}} 2^n (\sum_{J\in\bf{D}} |J| )
$$
where this time $\bf{D}$ ranges over those collections of disjoint dyadic intervals $J$ having 
the property that

$$\frac{1}{|J|}
\left\|\left(
\sum_{J'\in\J; J'\subseteq J}
\frac{|a_{J',k_0}^{(i)}|^2}{|J'|}\chi_{J'}(x)
\right)^{1/2}\right\|_{1,\infty}\geq 2^n.$$
\end{definition}

The following John - Nirenberg type inequality holds in this context, see \cite{cct}.

\begin{lemma}\label{l1}
Let $\J$ be a finite family of dyadic intervals as before. Then, for $i\neq j$ one has

$$\size^j_{i,\J}((a_J^{(i)})_J)\sim
\sup_{J\in\J}
\frac{1}{|J|^{1/2}}
(\sum_{J'\subseteq J} |a_{J'}^{(i)}|^2 )^{1/2}
$$
and similarly, if $j\neq 3$ one also has

$$\size^j_{3,k_0,\J}((a_{J,k_0}^{(3)})_J)\sim
\sup_{J\in\J}
\frac{1}{|J|^{1/2}}
(\sum_{J'\subseteq J} |a_{J',k_0}^{(3)}|^2 )^{1/2}.
$$
\end{lemma}

The following lemma, which has been proven in \cite{cct}, will also be very useful.

\begin{lemma}\label{l2}
Let $\J$ be as before and $i\neq j$. Then, for every $J\in\J$ one has the inequality

$$
\left\|\left(
\sum_{J'\in\J; J'\subseteq J}
\frac{|\langle f, \Phi^i_{J'}\rangle|^2}{|J'|}\chi_{J'}(x)
\right)^{1/2}\right\|_{1,\infty}\lesssim
\|f\widetilde{\chi}_J^N\|_1
$$
for every positive integer $N$, with the implicit constants depending on it.

\end{lemma}
 The following general inequality will play a fundamental role in our further estimates. It is
an abstract variant of the corresponding Proposition 3.6 in \cite{mptt:biparameter}.

\begin{proposition}\label{prop}
Let $\J$ be as before and $k_0\geq \#$. Then,

\begin{equation}\label{a}
|\Lambda_1(f_1, f_2, f_3, f_4)|\lesssim
\prod_{i=1}^3
(\size^j_{i,\J}((a_J^{(i)})_J))^{1-\theta_i}
(\energy^j_{i,\J}((a_J^{(i)})_J))^{\theta_i}
\end{equation}
for any $0\leq\theta_1, \theta_2,\theta_3<1$ such that $\theta_1+\theta_2+\theta_3 = 1$,
with the implicit constants depending on $\theta_i$ for $i=1,2,3$. Similarly, one also has

\begin{equation}\label{b}
|\Lambda_{1,k_0}(f_1, f_2, f_3, f_4)|\lesssim
\end{equation}

$$
(\size^j_{1,\J}((a_J^{(1)})_J))^{1-\theta_1}
(\size^j_{2,\J}((a_J^{(2)})_J))^{1-\theta_2}
(\size^j_{3,k_0\J}((a_J^{(3)})_J))^{1-\theta_3}\cdot
$$

$$
(\energy^j_{1,\J}((a_J^{(1)})_J))^{\theta_1}
(\energy^j_{2,\J}((a_J^{(2)})_J))^{\theta_2}
(\energy^j_{3,k_0\J}((a_J^{(3)})_J))^{\theta_3},
$$
for any $\theta_1, \theta_2, \theta_3$ exactly as before.
\end{proposition}

The proof of this Proposition will be presented later on. In the meantime we will take advantage
of it. In order to make it effective we would need to further estimate all these 
{\it sizes} and {\it energies} in terms of certain norms involving our functions
$f_1$, $f_2$, $f_3$, $f_4$.
The following lemma is an easy consequence of the previous definitions and of Lemma \ref{l2}
(see \cite{cct}).

\begin{lemma}\label{l3}
Let $E\subseteq \R$ be a set of finite measure, $i\neq 3$ and $f_{i+1}\in X(E)$. Then,

$$
\size^j_{i,\J}((a_J^{(i)})_J)\lesssim
\sup_{J\in\J}
\frac{1}{|J|}
\int_E \widetilde{\chi}_J^N dx,
$$
for every integer $N$, with the implicit constants depending on it.
\end{lemma}

Similarly, one also has

\begin{lemma}\label{l4}
With the same notations as in the previous lemma, we also have

$$\energy^j_{i,\J}((a_J^{(i)})_J)\lesssim |E|.$$
\end{lemma}

\begin{proof}
Let $n\in\Z$ and $\bf{D}$ be so that the suppremum in Definition \ref{se1} is attained. 
We also assume that $i\neq j$ (the case $i=j$ is in fact easier and is left to the reader).
Then, since the intervals $J\in\bf{D}$ are all disjoint, we can estimate the left hand side of
our inequality by

$$2^n(\sum_{J\in\bf{D}} |J|) =
2^n \|\sum_{J\in\bf{D}} \chi_J\|_1 =
2^n \|\sum_{J\in\bf{D}} \chi_J\|_{1,\infty}\lesssim
$$

$$\left\|
\sum_{J\in\bf{D}}
\frac{1}{|J|}
\|(
\sum_{J'\in\J; J'\subseteq J}
\frac{|a_{J'}^{(i)}|^2}{|J'|}\chi_{J'}(x))^{1/2}\|_{1,\infty}
\chi_J
\right\|_{1,\infty}\lesssim$$

$$\|
\sum_{J\in\bf{D}}
(\frac{1}{|J|}
\int |f_{i+1}|\widetilde{\chi}_J dx)
\chi_J
\|_{1,\infty}\lesssim$$

$$\|
\sum_{J\in\bf{D}}
(\frac{1}{|J|}
\int \chi_E \widetilde{\chi}_J dx)
\chi_J
\|_{1,\infty}\lesssim 
\|M(\chi_E)\|_{1,\infty}\lesssim |E|,$$
where $M$ is the Hardy - Littlewood maximal function and we also used Lemma \ref{l3}.
\end{proof}

We will also need

\begin{lemma}\label{l5}
Let $E_1, E_4\subseteq \R$ be sets of finite measure, $f_3\in X(E_3)$ and
$f_4\in X(E_4)$. Then,

$$\size^j_{3,\J}((a_J^{(3)})_J),\,\, 
\size^j_{3,k_0\J}((a_J^{(3)})_J)\lesssim
$$

$$\left(
\sup_{J\in\J}
\frac{1}{|J|}
\int_{E_1}\widetilde{\chi}_J^N dx \right)^{1-\theta}
\left(
\sup_{J\in\J}
\frac{1}{|J|}
\int_{E_4}\widetilde{\chi}_J^N dx \right)^{\theta},
$$
for any $0 <\theta < 1$ and for every positive integer $N$, with the implicit constants 
depending on them.
\end{lemma}

Similarly, we also have

\begin{lemma}\label{l6}
With the same notations as in the previous Lemma \ref{l5}, we have

$$\energy^j_{3,\J}((a_J^{(3)})_J),\,\, 
\energy^j_{3,k_0\J}((a_J^{(3)})_J)\lesssim
$$

$$\left(
\sup_{I\in\I_1}
\frac{1}{|I|}
\int_{E_1}\widetilde{\chi}_I^N dx \right)^{1-\theta_1}
\left(
\sup_{I\in\I_1}
\frac{1}{|I|}
\int_{E_4}\widetilde{\chi}_I^N dx \right)^{1-\theta_2}
|E_1|^{\theta_1} |E_4|^{\theta_2},
$$
for any $0\leq \theta_1, \theta_2 < 1$ with $\theta_1 + \theta_2 = 1$ and for every integer
$N$, with the implicit constants depending on them.
\end{lemma}

The proofs of these two lemmas will be presented later on. In the meantime, we will take 
advantage
of them, in order to complete the proof of our Theorem \ref{modelthm}.

\section{Estimates for $T_1$ and $T_{1,k_0}$ near $A_4$}

In this section we start the proof of Theorem \ref{modelthm}. Clearly, due to symmetry
considerations, it is enough to analyze the case of $T_1$ and $T_{1, k_0}$, the case
of $T_1$ and $T_{2, k_0}$ being similar.

Let now $(p_1, p_2, p_3, p_4)$ be so that $(1/p_1, 1/p_2, 1/p_3, 1/p_4)\in \bf{D}$ and is
arbitrarily close to $A_4$ which has coordinates $(1,1,1,-2)$. Let also
$E_1$, $E_2$, $E_3$, $E_4\subseteq \R$ be measurable sets of finite measure. By
scaling invariance, we can also assume that $|E_4| = 1$. Our goal is to construct
a subset $E'_4\subseteq E_4$ with $|E'_4| \sim 1$ and so that

\begin{equation}\label{16}
|\Lambda_1(f_1, f_2, f_3, f_4)|, \,\,
|\Lambda_{1,k_0}(f_1, f_2, f_3, f_4)|\lesssim
|E_1|^{1/p_1}
|E_2|^{1/p_2}
|E_3|^{1/p_3}
\end{equation}
for every $f_i\in X(E_i)$, $i=1,2,3$ and $f_4\in X(E'_4)$. As in \cite{mptt:biparameter}
define first the exceptional set $\Omega$ by

$$\Omega:= \bigcup_{j=1}^3
\{ M(\frac{\chi_{E_j}}{|E_j|}) > C \}
$$
and observe that $|\Omega| < 1/2$ if $C$ is a big enough constant. Then, set 
$E'_4:= E_4\setminus \Omega$ which clearly has the property that $|E'_4|\sim 1$.

Now, we decompose the sets $\J_1$ and $\I_1$ as

$$\J_1:= \bigcup_{d\geq 0} \J_1^d$$

$$\I_1:= \bigcup_{d'\geq 0} \I_1^{d'}$$
where $\J_1^d$ is the set of all intervals $J\in\J_1$ with the property that

$$\left(1 + \frac{\dist (J,\Omega^c)}{|J|}\right) \sim 2^d$$
and $\I_1^{d'}$ is the set of all intervals $I\in\I_1$ with the property that

$$\left(1 + \frac{\dist (I,\Omega^c)}{|J|}\right) \sim 2^{d'}.$$

From the definition of $\Omega$ we have

\begin{equation}\label{17}
\frac{1}{|J|}
\int_{E_j}\widetilde{\chi}_J dx
\lesssim 2^d |E_j|
\end{equation}
for $j=1,2,3$ and since obviously the left hand side of (\ref{17}) is also smaller than
$1$, it follows that

$$\frac{1}{|J|}
\int_{E_j}\widetilde{\chi}_J dx
\lesssim 2^{\alpha d} |E_j|^{\alpha}
$$
for every $0\leq\alpha\leq 1$ and $j=1,2,3$.

Similarly, we also have that

$$\frac{1}{|I|}
\int_{E_j}\widetilde{\chi}_I dx
\lesssim 2^{\beta d'} |E_j|^{\alpha}
$$
for every $0\leq\beta\leq 1$ and $j=1,2,3$.

On the other hand, since $E'_4\subseteq\Omega^c$ we also know that

$$\frac{1}{|J|}
\int_{E'_4}\widetilde{\chi}_J dx
\lesssim 2^{-Nd}
$$
and

$$\frac{1}{|I|}
\int_{E'_4}\widetilde{\chi}_I dx
\lesssim 2^{-Nd'}
$$
for any integer $N > 0$. Using now our previous lemmas together with all these observations,
we obtain the estimates

$$
\size^j_{1, \J_1^d}((a_J^{(1)})_J)\lesssim 2^{d\alpha_2} |E_2|^{\alpha_2}
$$

$$
\size^j_{2, \J_1^d}((a_J^{(2)})_J)\lesssim 2^{d\alpha_3} |E_3|^{\alpha_3}
$$

$$
\size^j_{3, \J_1^d}((a_J^{(3)})_J), \,\,
\size^j_{3,k_0 \J_1^d}((a_{J,k_0}^{(3)})_J)\lesssim
(2^{d\alpha_1} |E_1|^{\alpha_1})^{1-\theta} (2^{-Nd})^{\theta}
$$
and similarly,

$$
\energy^j_{1, \J_1^d}((a_J^{(1)})_J)\lesssim |E_2|
$$

$$
\energy^j_{2, \J_1^d}((a_J^{(2)})_J)\lesssim |E_3|
$$

$$
\energy^j_{3, \J_1^d}((a_J^{(3)})_J), \,\,
\energy^j_{3,k_0 \J_1^d}((a_{J,k_0}^{(3)})_J)\lesssim
(2^{d'\beta_1} |E_1|^{\beta_1})^{1-\theta'_1} (2^{-Nd'})^{1-\theta'_2}|E_1|^{\theta'_1}
$$
whenever $0\leq\alpha_1, \alpha_2, \alpha_3, \beta_1\leq 1$, $0 < \theta < 1$ and
$0\leq \theta'_1, \theta'_2 < 1$ with $\theta'_1 + \theta'_2 = 1$.

By using now Proposition \ref{prop} we deduce that for any 
$0\leq\theta_1, \theta_2, \theta_3 < 1$ with $\theta_1 + \theta_2 +\theta_3 = 1$, 
one can estimate
the left hand side of (\ref{16}) by

$$
(2^{d\alpha_2}|E_2|^{\alpha_2})^{1-\theta_1}
(2^{d\alpha_3}|E_3|^{\alpha_3})^{1-\theta_2}
[(2^d|E_1|)^{1-\theta} (2^{-Nd})^{\theta}]^{1-\theta_3}\cdot
$$

$$|E_2|^{\theta_1} |E_3|^{\theta_2}
[(2^{d'}|E_1|)^{1-\theta'_1} (2^{-Nd'})^{1-\theta'_2} |E_1|^{\theta'_1}]^{\theta_3} =
$$

$$|E_1|^{(1-\theta)(1-\theta_3) + \theta_3}\cdot
|E_2|^{\alpha_2(1-\theta_1) + \theta_1}\cdot
|E_3|^{\alpha_3(1-\theta_2) + \theta_2}\cdot
2^{-u d}\cdot 2^{-v d'}
$$
where $u, v$ are both positive numbers depending on all these parameters and also on $N$.

Now, if one takes $\theta_1$ very close to $0$ and $\alpha_2, \alpha_3$ very close to $1$,
one can then define $1/p_1:= (1-\theta)(1-\theta_3) + \theta_3$,
$1/p_2:= \alpha_2(1-\theta_1) + \theta_1$ and $1/p_3:= \alpha_3(1-\theta_2) + \theta_2$
and they can be chosen as close as we want to the point $(1,1,1)$.

In the end, one can sum over $d, d'\geq 0$ if our constant $N$ is big enough.

A similar argument proves the desired estimates for $T_1^{*1}$ and
$T_{1, k_0}^{*1}$ near the points $A_{11}$ and $A_{12}$.

\section{Estimates for $T_1^{*3}$ and $T_{1,k_0}^{*3}$ near $A_{31}$ and $A_{32}$ }

The proof uses similar ideas as in the argument in the previous section.

Let $(p_1, p_2, p_3, p_4)$ so that $(1/p_1, 1/p_2, 1/p_3, 1/p_4)\in \bf{D}$ and is
arbitrarily close to either $A_{31}$ or $A_{32}$. Consider $E_1$, $E_2$, $E_3$, 
$E_4\subseteq \R$ measurable sets of finite measure and assume as before that $|E_3| = 1$. 
Our task
is to construct a subset $E'_3\subseteq E_3$ with $|E'_3|\sim 1$ so that

\begin{equation}\label{18}
|\Lambda_1(f_1, f_2, f_3, f_4)|, \,\,
|\Lambda_{1,k_0}(f_1, f_2, f_3, f_4)|\lesssim
|E_1|^{1/p_1}
|E_2|^{1/p_2}
|E_4|^{1/p_4}
\end{equation}
for every $f_i\in X(E_i)$, $i=1,2,4$ and $f_3\in X(E'_3)$.

Define the exceptional set 

$$\Omega:= \{ M(\frac{\chi_{E_1}}{|E_1|}) > C \} \cup
\{ M(\frac{\chi_{E_2}}{|E_2|}) > C \} \cup
\{ M(\frac{\chi_{E_4}}{|E_4|}) > C \}
$$
and then set $E'_3:= E_3\setminus \Omega$ for a sufficiently large constant $C > 0$.

With the same notations as in Section 6, we obtain the estimates (this time there is no need
to decompose $\I_1$ as there)

$$
\size^j_{1, \J_1^d}((a_J^{(1)})_J)\lesssim 2^{d\alpha_2} |E_2|^{\alpha_2}
$$

$$
\size^j_{2, \J_1^d}((a_J^{(2)})_J)\lesssim 2^{-Nd}
$$

$$
\size^j_{3, \J_1^d}((a_J^{(3)})_J), \,\,
\size^j_{3,k_0 \J_1^d}((a_{J,k_0}^{(3)})_J)\lesssim
(2^{d\alpha_1} |E_1|^{\alpha_1})^{1-\theta} (2^{d\alpha_4} |E_4|^{\alpha_4} )^{\theta}
$$
and similarly,

$$
\energy^j_{1, \J_1^d}((a_J^{(1)})_J)\lesssim |E_2|
$$

$$
\energy^j_{2, \J_1^d}((a_J^{(2)})_J)\lesssim 1
$$

$$
\energy^j_{3, \J_1^d}((a_J^{(3)})_J), \,\,
\energy^j_{3,k_0 \J_1^d}((a_{J,k_0}^{(3)})_J)\lesssim
|E_1|^{\widetilde{\theta}_1}|E_4|^{\widetilde{\theta}_2}
$$
whenever $0\leq \alpha_1, \alpha_2, \alpha_4\leq 1$, $0 <\theta < 1$ and 
$0\leq \widetilde{\theta}_1, \widetilde{\theta}_2 < 1$ with
$\widetilde{\theta}_1 + \widetilde{\theta}_2  = 1$.

Then, by applying Proposition \ref{prop} we obtain the following estimates for the left hand side
of (\ref{18})

$$
(2^{d\alpha_2}|E_2|^{\alpha_2})^{1-\theta_1}
(2^{-Nd})^{1-\theta_2}
[(2^d|E_1|)^{1-\theta} (2^{d}|E_4|)^{\theta}]^{1-\theta_3}\cdot
$$

$$|E_2|^{\theta_1}
[|E_1|^{\widetilde{\theta}_1}|E_4|^{\widetilde{\theta}_2} ]^{\theta_3} =
$$

$$|E_1|^{(1-\theta)(1-\theta_3) + \widetilde{\theta}_1\theta_3}\cdot
|E_2|^{\alpha_2(1-\theta_1) + \theta_1}\cdot
|E_4|^{\theta(1-\theta_3) + \widetilde{\theta}_2\theta_3}\cdot
2^{-u d}
$$
where again $u$ is a positive number depending on all these parameters.

Then, we define $1/p_1:= (1-\theta)(1-\theta_3) + \widetilde{\theta}_1\theta_3$,
$1/p_2:= \alpha_2(1-\theta_1) + \theta_1$ and 
$1/p_4:= \theta(1-\theta_3) + \widetilde{\theta}_2\theta_3$ and since
$(1-\theta)(1-\theta_3) + \widetilde{\theta}_1\theta_3 +
\theta(1-\theta_3) + \widetilde{\theta}_2\theta_3 = 1$, one can easily check
that $p_2$ can be chosen very close to $1$ (by chosing $\alpha_2$ close to $1$)
and the pair $(1/p_1, 1/p_4)$ very close either to $(0,1)$ or $(1,0)$ which is what we
wanted. And in the end we sum over $d\geq 0$ since $u$ remains positive if we chose $N$ big
enough.

A similar argument proves the required estimates for the operators
$T_1^{*2}$ and $T^{*2}_{1,k_0}$ near $A_{21}$ and $A_{22}$.

\section{Proof of Proposition \ref{prop}}

This section is devoted to the proof of Proposition \ref{prop}. As we pointed out earlier,
this proposition is a more abstract version of the corresponding Proposition 3.6 in
\cite{mptt:biparameter}. Its proof is similar and we include it here for completeness and also 
for the reader's convenience.

\begin{proposition}\label{prop1}
Let $\J$ be a finite family of dyadic intervals, $\J'$ a subset of $\J$, $i=1,2,3$, $n_0\in\Z$
and assume that

$$
\size^j_{i,\J'}((a_J^{(i)})_J)\leq 2^{-n_0}\energy^j_{i,\J}((a_J^{(i)})_J).
$$
Then, there exists a decompostion $\J' = \J''\cup \J'''$ such that

\begin{equation}\label{19}
\size^j_{i,\J''}((a_J^{(i)})_J)\leq 2^{-n_0-1}\energy^j_{i,\J}((a_J^{(i)})_J)
\end{equation}
and so that $\J'''$ can be written as a disjoint union of subsets $T\in\T$ such that
for every $T\in\T$ there exists a dyadic interval $J_T\in\J$ having the property that every
$J\in T$ satisfies $J\subseteq J_T$ and also such that

\begin{equation}\label{20}
\sum_{T\in\T} |J_T| \lesssim 2^{n_0}.
\end{equation}
\end{proposition}

\begin{proof}
\underline{Case $1$}: $i=j$. First, chose an interval $J\in\J'$ having the property
that $|J|$ is as big as possible and so that

\begin{equation}\label{21}
\frac{|a_J^{(i)}|}{|J|} > 2^{-n_0-1}\energy^j_{i,\J}((a_J^{(i)})_J).
\end{equation}
Then, collect all the intervals $J'\in\J'$ with $J'\subseteq J$ into a set called $T$. After
this, define $J_T:= J$ and look at the remaining intervals in $\J'\setminus T$ and repeat the
procedure. Since there are finitely many such dyadic intervals, the procedure ends after finitely
many steps producing the subsets $T\in\T$. Define $\J''':= \cup_{T\in\T} T$
and $\J'':= \J\setminus\J'''$. Now clearly, by construction, the inequality (\ref{19})
is satisfied and it only remains to check (\ref{20}).

Since the intervals $(J_T)_{T\in\T}$ are all disjoint by construction, we deduce from
(\ref{21}) and Definition \ref{se1} that

$$2^{-n_0}\energy^j_{i,\J}((a_J^{(i)})_J) (\sum_{T\in\T} |J_T| )\lesssim
\energy^j_{i,\J}((a_J^{(i)})_J)
$$
which is equivalent to our desired estimate (\ref{20}).

\underline{Case $2$}: $i\neq j$. The procedure of selecting the intervals is very similar. 
The only difference is that this time, we pick intervals $J'\in\J'$ so that $|J|$ is again
as big as possible, but having the property that

$$
\frac{1}{|J|}
\left\|\left(
\sum_{J'\in\J; J'\subseteq J}
\frac{|a_{J'}^{(i)}|^2}{|J'|}\chi_{J'}(x)\right)^{1/2}\right\|_{1,\infty}
 > 2^{-n_0-1}\energy^j_{i,\J}((a_J^{(i)})_J).
$$
After this the argument is identical to the one we described before.
\end{proof}

Similarly, we also have

\begin{proposition}\label{prop2}
Using the same notations as in the previous Proposition \ref{prop1}, assume that

$$
\size^j_{3, k_0,\J'}((a_J^{(3)})_J)\leq 2^{-n_0}\energy^j_{3, k_0,\J}((a_J^{(3)})_J).
$$
Then, there exists a decompostion $\J' = \J''\cup \J'''$ as before, such that

\begin{equation}
\size^j_{3, k_0,\J''}((a_J^{(3)})_J)\leq 2^{-n_0-1}\energy^j_{3,k_0,\J}((a_J^{(3)})_J)
\end{equation}
and so that $\J'''$ can be written as a disjoint union of subsets $T\in\T$ such that
for every $T\in\T$ there exists a dyadic interval $J_T\in\J$ having the property that every
$J\in T$ satisfies $J\subseteq J_T$ and also such that

\begin{equation}
\sum_{T\in\T} |J_T| \lesssim 2^{n_0}.
\end{equation}
\end{proposition}

By iterating these two propositions, we obtain the following corollaries.

\begin{corollary}\label{cor1}
Let $i=1,2,3$ and $\J$ be a finite family of dyadic intervals. Then, there exists a partition

$$\J = \bigcup_{n\in\Z}\J^{n,i}$$
such that for every $n\in\Z$ we have

$$\size^j_{i,\J^{n,i}}((a_J^{(i)})_J)\leq min (2^{-n}\energy^j_{i,\J}((a_J^{(i)})_J),
\size^j_{i,\J}((a_J^{(i)})_J)).
$$
Also, we can write each $\J^{n,i}$ as a disjoint union of subsets $T\in\T_n^i$ as before,
having the property that

$$\sum_{T\in\T_n^i} |J_T|\lesssim 2^n.$$
\end{corollary}

\begin{corollary}\label{cor2}
Let $\J$ be a finite family of dyadic intervals. Then, there exists a partition

$$\J = \bigcup_{n\in\Z}\J^{n}$$
such that for every $n\in\Z$ we have

$$\size^j_{3. k_0,\J^{n}}((a_J^{(3)})_J)\leq min (2^{-n}\energy^j_{3, k_0,\J}((a_J^{(3)})_J),
\size^j_{3, k_0,\J}((a_J^{(3)})_J)).
$$
Also, we can write each $\J^{n}$ as a disjoint union of subsets $T\in\T_n$ as before,
having the property that

$$\sum_{T\in\T_n} |J_T|\lesssim 2^n.$$
\end{corollary}

Having all these decompositions available, we can now start the actual proof of Proposition
\ref{prop}. We will only present the proof of the first inequalty (\ref{a}), the proof
of (\ref{b}) being similar.

As in \cite{mptt:biparameter}, since $j$ is fixed anyways, we will write for simplicity 
$S_i:= \size^j_{i,\J}((a_J^{(i)})_J)$ and 
$E_i:= \energy^j_{i,\J}((a_J^{(i)})_J) $, for
$i=1,2,3$. If we apply Corollary \ref{cor1} to our collection $\J$, we obtain
a decomposition
 
$$\J=\bigcup_n \J^{n,i}$$
such that each $\J^{n,i}$ can be written as a union of subsets in $\T^i_n$ with the properties 
described in 
Corollary \ref{cor1}.
Consequently, one can estimate the left hand side of our inequality (\ref{a}) as

\begin{equation}\label{in1}
\sum_{n_1, n_2, n_3} \sum_{T\in \T^{n_1, n_2, n_3}} 
\sum_{J\in T}\frac{1}{|J|^{1/2}} |a_J^{(1)}||a_J^{(2)}||a_J^{(3)}|
\end{equation}
where $\T^{n_1, n_2, n_3}:= \T^1_{n_1} \cap \T^2_{n_2} \cap \T^3_{n_3}$.  

Fix such a $T$ and look at the corresponding inner term in (\ref{in1}). It can be estimated by

$$\sup_{J\in T}\frac{|a_J^{(j)}|}{|J|^{1/2}}
\prod_{i\neq j}(\sum_{J\in T}|a_J^{(i)}|^2 )^{1/2} =
$$

$$\sup_{J\in T}\frac{|a_J^{(j)}|}{|J|^{1/2}}
\left(\prod_{i\neq j}\frac{1}{|J_T|^{1/2}}(\sum_{J\in T}|a_J^{(i)}|^2 )^{1/2}\right) |J_T|\lesssim
$$

$$\left(\prod_{i=1}^3\size^j_{i, T}((a_J^{(i)})_J)\right) |J_T|,$$
by also using the John-Nirenberg inequality in Lemma \ref{l1}. 

In particular, we can estimate (\ref{in1}) further by

\begin{equation}\label{in2}
E_1 E_2 E_3\sum_{n_1, n_2, n_3} 2^{-n_1} 2^{-n_2} 2^{-n_3} \sum_{T\in\T^{n_1, n_2, n_3}}|I_T|
\end{equation}
where, according to the same Corollary \ref{cor1}, the summation goes over those 
$n_1, n_2, n_3 \in \Z$ having the property that

\begin{equation}\label{must}
2^{-n_j}\lesssim \frac{S_j}{E_j}.
\end{equation}
On the other hand, Corollary \ref{cor1} allows us to estimate the inner sum in 
(\ref{in2}) in three different ways,
namely

$$\sum_{T\in\T^{n_1, n_2, n_3}} |I_T|\lesssim 2^{n_1}, 2^{n_2}, 2^{n_3}$$
and so, as a consequence, we can also write

\begin{equation}\label{in3}
\sum_{T\in\T^{n_1, n_2, n_3}} |I_T|\lesssim 2^{n_1\theta_1} 2^{n_2\theta_2} 2^{n_3\theta_3}
\end{equation}
whenever $0\leq \theta_1, \theta_2, \theta_3 <1$ with $\theta_1+ \theta_2+ \theta_3= 1$. 
Using (\ref{in3}) and (\ref{must}),
one can estimate (\ref{in2}) again by

$$E_1 E_2 E_3 \sum_{n_1, n_2, n_3}2^{-n_1(1-\theta_1)} 2^{-n_2(1-\theta_2)} 2^{-n_3(1-\theta_3)}\lesssim$$

$$E_1 E_2 E_3 
(\frac{S_1}{E_1})^{1-\theta_1}
(\frac{S_2}{E_2})^{1-\theta_2} 
(\frac{S_2}{E_2})^{1-\theta_3}= \prod_{j=1}^3 S_j^{1-\theta_j} \prod_{j=1}^3 E_j^{\theta_j},$$
and this ends the proof.

\section{Proof of Lemma \ref{l5}}

\underline{Case $I$}: Estimates for $\size^j_{3,\J}((a^{(3)})_J)$.

These are essentially known (see \cite{mtt:fourierbiest}). We include a slightly different 
proof here for 
completeness and also since the same argument has enough flexibility to also handle the case of
$\size^j_{3,k_0,\J}((a^{(3)})_J)$ later on. There are two subcases.

\underline{Case $I_1$}: $j\neq 3$.

Fix $J_0\in\J$. Clearly, to prove our estimates it is enough to show that

\begin{equation}\label{22}
\left\|\left(
\sum_{J\subseteq J_0}
\frac{|a^{(3)}_J|^2}{|J|}\chi_J(x) \right)^{1/2}\right\|_{1,\infty}
\lesssim \|f_3 \widetilde{\chi}_{J_0}^N\|_p\cdot\|f_4 \widetilde{\chi}_{J_0}^N\|_q
\end{equation}
whenever $1 < p, q <\infty$ with $1/p + 1/q = 1$. Let us now recall that $a^{(3)}_J$ is
defined by

\begin{equation}\label{23}
a_J^{(3)}:=
\langle
\sum_{I\in\I_1; \omega_J^3\cap \omega_I^2\neq\emptyset; |\omega_J^3|\leq |\omega_I^2|}
\frac{1}{|I|^{1/2}}
\langle f_1, \Phi_I^1\rangle
\langle f_4, \Phi_I^3\rangle
\Phi_I^2, \Phi_J^3\rangle.
\end{equation}
Define the collection $\widetilde{\I}$ to be the set of all dyadic intervals $I\in\I_1$ for
which there exists $J\in\J$ with the property that 
$\omega_J^3\cap \omega_I^2\neq\emptyset$ and $|\omega_J^3|\leq |\omega_I^2|$. We claim that

\begin{equation}\label{24}
a^{(3)}_J = \langle B(f_1, f_4), \phi_J^3\rangle
\end{equation}
where $B(f_1, f_4)$ is defined by

\begin{equation}\label{25}
B(f_1, f_4):= 
\sum_{I\in\widetilde{\I}}
\frac{1}{|I|^{1/2}}
\langle f_1, \Phi_I^1\rangle
\langle f_4, \Phi_I^3\rangle
\Phi_I^2 (x).
\end{equation}
To check the claim, let us observe that for each $I\in\widetilde{\I}$

$$\langle\Phi^2_I,\Phi^3_J\rangle \neq 0\,\,\,\text{iff}\,\,\,
\omega^3_J\cap\omega^2_I\neq\emptyset.$$
There are two possibilities: either $|\omega^3_J|\leq |\omega^2_I|$ which is acceptable
by (\ref{23}), or $|\omega^2_I| < |\omega^3_J|$. We then make the claim that this last
situation cannot occur. Indeed, since $\omega^2_I$ is symmetric with respect to the origin,
that would imply that $0\in 3 \omega^3_J$ which is clearly false, since by 
Definition \ref{lacunary} one has $0\notin 5 \omega^3_J$.

Using now (\ref{24}) together with Lemma \ref{l2} it follows that to prove (\ref{22})
it is enough to prove that

\begin{equation}\label{26}
\|B(f_1, f_4)\widetilde{\chi}_{J_0}^{N'}\|_1
\lesssim
\|f_3 \widetilde{\chi}_{J_0}^N\|_p\cdot\|f_4 \widetilde{\chi}_{J_0}^N\|_q
\end{equation}
By scaling invariance, we may assume without loss of generality that $|J_0| = 1$.
Our plan is to prove a slightly weaker version of (\ref{26}), namely to prove that

\begin{equation}\label{27}
\|B(f_1, f_4)\chi_J\|_1
\lesssim
\|f_3 \widetilde{\chi}_{J}^N\|_p\cdot\|f_4 \widetilde{\chi}_{J}^N\|_q
\end{equation}
for every dyadic interval $J\subseteq \R$ of length $1$. We now prove that if we assume
(\ref{27}) then (\ref{26}) follows quite easily.

To see this, consider a partition of the real line with disjoint intervals of length $1$
$(J_n)_{n\in\Z^*}$ so that

$$(\bigcup_{n\in\Z^*}J_n)\cup J_0 = \R.$$
Then, estimate the left hand side of (\ref{26}) by

$$\|B(f_1, f_4)\widetilde{\chi}_{J_0}^{N'}\|_1\lesssim
\|B(f_1, f_4)\chi_{J_0}\|_1 + 
\sum_{n\in\Z^*}
\|B(f_1, f_4)\widetilde{\chi}_{J_0}^{N'}\chi_{J_n}\|_1\lesssim $$

$$\|B(f_1, f_4)\chi_{J_0}\|_1 + \sum_{n\in\Z^*}\frac{1}{|n|^{N'}}
\|B(f_1, f_4)\chi_{J_n}\|_1.$$
The first term clearly satisfies the desired estimates. The second one can be further majorized
using (\ref{27}) by

$$\sum_{n\in\Z^*}\frac{1}{|n|^{N'}}
\|f_1\widetilde{\chi}_{J_n}^{N''}\|_p\cdot
\|f_4\widetilde{\chi}_{J_n}^{N''}\|_q\lesssim
$$

$$\left(\sum_{n\in\Z^*}\frac{1}{|n|^{N'}}
\|f_1\widetilde{\chi}_{J_n}^{N''}\|_p^p\right)^{1/p}\cdot
\left(\sum_{n\in\Z^*}\frac{1}{|n|^{N'}}
\|f_4\widetilde{\chi}_{J_n}^{N''}\|_q^q\right)^{1/q}\lesssim
$$

$$\left(\int_{\R} |f_1|^p(\sum_{n\in\Z^*}\frac{1}{|n|^{N'}}\widetilde{\chi}_{J_n}^{pN''}) dx
\right)^{1/p}\cdot
\left(\int_{\R} |f_4|^q(\sum_{n\in\Z^*}\frac{1}{|n|^{N'}}\widetilde{\chi}_{J_n}^{qN''}) dx
\right)^{1/q}\lesssim
$$

$$\|f_1\widetilde{\chi}^N_{J_0}\|_p\cdot
\|f_1\widetilde{\chi}^N_{J_0}\|_q
$$
if $N'$ is big enough. It remains to prove (\ref{27}).

\underline{Case $I_{1a}$}: $\supp f_1$, $\supp f_4 \subseteq 5J$.

In this case, our inequality (\ref{27}) follows from the known estimates on discrete paraproducts
(see for instance \cite{mptt:biparameter} ).

\underline{Case $I_{1b}$}: Either $\supp f_1\subseteq (5J)^c$ or $\supp f_4\subseteq (5J)^c$.

Assume for instance that $\supp f_1\subseteq (5J)^c$. Then, we decompose $B(f_1, f_4)$ as

$$B(f_1, f_4) = B'(f_1, f_4) + B''(f_1, f_4)$$
where

$$B'(f_1, f_4):= \sum_{I\in\widetilde{\I}; I\cap 5J\neq\emptyset}
\frac{1}{|I|^{1/2}}
\langle f_1, \Phi_I^1\rangle
\langle f_4, \Phi_I^3\rangle
\Phi_I^2 $$
and

$$B''(f_1, f_4):= \sum_{I\in\widetilde{\I}; I\cap 5J = \emptyset}
\frac{1}{|I|^{1/2}}
\langle f_1, \Phi_I^1\rangle
\langle f_4, \Phi_I^3\rangle
\Phi_I^2 $$
By our reduction ( $|J_0| = 1$) we observe that the lengths of our intervals $I$ are all
smaller than $1$.

For $h\in L^{\infty}$, $\|h\|_{\infty}\leq 1$ one can write

$$|\int_{\R} B'(f_1, f_4)(x) h(x) \chi_J(x) dx |\lesssim $$

$$\sum_{I\in\widetilde{\I}; I\cap 5J\neq\emptyset}
\frac{1}{|I|^{1/2}}
|\langle f_1, \Phi_I^1\rangle|
|\langle f_4, \Phi_I^3\rangle|
|\langle h\chi_J,\Phi_I^2\rangle|= $$

$$\sum_{k=0}^{\infty}
\sum_{I\in\widetilde{\I}; I\cap 5J\neq\emptyset; |I| = 2^{-k}}
2^{k/2}
|\langle f_1, \Phi_I^1\rangle|
|\langle f_4, \Phi_I^3\rangle|
|\langle h\chi_J,\Phi_I^2\rangle|= $$

\begin{equation}\label{28}
\sum_{k=0}^{\infty}
\sum_{I\in\widetilde{\I}; I\cap 5J\neq\emptyset; |I| = 2^{-k}}
2^{2k}
|\langle f_1, 2^{-k/2}\Phi_I^1\rangle|
|\langle f_4, 2^{-k/2}\Phi_I^3\rangle|
|\langle h\chi_J, 2^{-k/2}\Phi_I^2\rangle|
\end{equation}
and observe that all the functions $2^{-k/2}\Phi_I^1$, $2^{-k/2}\Phi_I^3$ and
$2^{-k/2}\Phi_I^2$ are $L^{\infty}$- normalized. Then, we estimate (\ref{28}) by

$$\sum_{k=0}^{\infty}2^{2k}
\left(\sup_{I\cap 5J\neq\emptyset; |I| = 2^{-k}}|\langle f_1, 2^{-k/2}\Phi_I^1\rangle|\right)
\left(\sup_{I\cap 5J\neq\emptyset; |I| = 2^{-k}}|\langle f_4, 2^{-k/2}\Phi_I^3\rangle|\right)
\cdot$$

$$|\langle h\chi_J, \sum_{I\in\widetilde{\I}; I\cap 5J\neq\emptyset; |I| = 2^{-k}}
\widetilde{\chi}_I^N \rangle|\lesssim
$$

$$\sum_{k=0}^{\infty} 2^{2k} 2^{-100k}
\|f_1\widetilde{\chi}^N_J\|_1\cdot
\|f_4\widetilde{\chi}^N_J\|_1\lesssim
\|f_1\widetilde{\chi}^N_J\|_p\cdot
\|f_4\widetilde{\chi}^N_J\|_q.
$$

Similarly, one can also write

$$|\int_{\R} B''(f_1, f_4)(x) h(x) \chi_J(x) dx |\lesssim $$

$$\sum_{I\in\widetilde{\I}; I\cap 5J = \emptyset}
\frac{1}{|I|^{1/2}}
|\langle f_1, \Phi_I^1\rangle|
|\langle f_4, \Phi_I^3\rangle|
|\langle h\chi_J,\Phi_I^2\rangle|= $$

$$
\sum_{k=0}^{\infty}
\sum_{I\in\widetilde{\I}; I\cap 5J = \emptyset; |I| = 2^{-k}}
2^{2k}
|\langle f_1, 2^{-k/2}\Phi_I^1\rangle|
|\langle f_4, 2^{-k/2}\Phi_I^3\rangle|
|\langle h\chi_J, 2^{-k/2}\Phi_I^2\rangle|\lesssim
$$

$$\sum_{k=0}^{\infty}
2^{2k}
\sum_{I\in\widetilde{\I}; I\cap 5J = \emptyset; |I| = 2^{-k}}
\dist(I, J)^{2N}
\|f_1\widetilde{\chi}^N_J\|_1\cdot
\|f_4\widetilde{\chi}^N_J\|_1\cdot
(\frac{\dist(I, J)}{|I|})^{-N'}\lesssim
$$

\begin{equation}\label{29}
\|f_1\widetilde{\chi}^N_J\|_p\cdot
\|f_4\widetilde{\chi}^N_J\|_q\cdot
\sum_{k=0}^{\infty}2^{-(N'-2)k}
\sum_{I\in\widetilde{\I}; I\cap 5J = \emptyset; |I| = 2^{-k}}
(\dist(I, J))^{-(N'-2N)}.
\end{equation}
Now, if $N'$ is much bigger than $2N$ then the inner sum in (\ref{29}) is smaller than

$$\sum_{n=0}^{\infty}
\frac{1}{(1+ n 2^{-k})^{N'-2N}} =
2^{k(N'-2N)}\sum_{n=0}^{\infty}
\frac{1}{(2^k+ n)^{N'-2N}}\lesssim
$$

$$2^{k(N'-2N)}\int_{2^k}^{\infty}\frac{1}{x^{N'-2N}} dx \lesssim 2^k$$
and this makes the geometric series in (\ref{29}) convergent. We are then left with
\underline{Case $I_2$} when $j=2$ but this clearly follows by the same arguments.

\underline{Case $II$}: Estimates for $\size^j_{3, k_0,\J}((a^{(3)}_J))$.

The argument follows the same ideas as before. There are several subcases.

\underline{Case $II_1$}: $j\neq 3$.

Fix as before $J_0\in\J$. Clearly, to prove our estimates it is enough to show that

\begin{equation}\label{30}
\left\|\left(
\sum_{J\subseteq J_0}
\frac{|a^{(3)}_{J,k_0}|^2}{|J|}\chi_J(x) \right)^{1/2}\right\|_{1,\infty}
\lesssim \|f_3 \widetilde{\chi}_{J_0}^N\|_p\cdot\|f_4 \widetilde{\chi}_{J_0}^N\|_q
\end{equation}
whenever $1 < p, q <\infty$ with $1/p + 1/q = 1$. Let us now recall that $a^{(3)}_{J,k_0}$ is
defined by the formula

\begin{equation}\label{31}
a_J^{(3)}:=
\langle
\sum_{I\in\I_1; \omega_J^3\cap \omega_I^2\neq\emptyset; 2^{k_0}|\omega_J^3|\sim |\omega_I^2|}
\frac{1}{|I|^{1/2}}
\langle f_1, \Phi_I^1\rangle
\langle f_4, \Phi_I^3\rangle
\Phi_I^2, \Phi_J^3\rangle.
\end{equation}
Since the frequency intervals $\omega_I^2$ and $\omega_J^3$ depend only on the scales
$|I|$ and $|J|$ respectively (see Section $3$) it follows  that by a certain refinement
we can assume that given $|J|$ there exists only one $|I|$ so that
$\omega^3_{|J|}\cap\omega^2_{|I|}\neq\emptyset$ and
$2^{k_0}|\omega_{|J|}^3|\sim |\omega_{|I|}^2|$. Fix now such a pair of dyadic intervals
$I$ and $J$. Then, by Plancherel, we have

\begin{equation}\label{31'}
\langle\Phi^2_I,\phi^3_J\rangle =
\langle\widehat{\Phi^2_I}, \widehat{\phi^3_J}\rangle.
\end{equation}
Since $|J|\sim 2^{k_0} |I|$, pick a Schwartz function $\Psi_{|I|, k_0}$ so that
$\supp\widehat{\Psi_{|I|, k_0}}\subseteq 2\omega^3_{|J|}$ and 
$\widehat{\Psi_{|I|, k_0}}\equiv 1$ on $\omega^3_{|J|}$.

Then, (\ref{31'}) equals

$$\langle\widehat{\Phi^2_I}, \widehat{\phi^3_J}\cdot\widehat{\Psi_{|I|, k_0}}
\rangle =
\langle \widehat{\Phi^2_I *\Psi_{|I|, k_0} }, \widehat{\Phi^3_J}\rangle =
$$

$$2^{-k_0/2}
\langle \widehat{2^{k_0/2}\Phi^2_I *\Psi_{|I|, k_0} }, \widehat{\Phi^3_J}\rangle =
2^{-k_0/2}\langle\widetilde{\Phi}^2_I, \Phi^3_J\rangle
$$
where

$$\widetilde{\Phi}^2_I:= 2^{k_0/2}\Phi^2_I *\Psi_{|I|, k_0}$$
and it is not difficult to observe that $\widetilde{\Phi}^2_I$ is an $L^2$- normalized
bump adapted to $\widetilde{I}$, where $\widetilde{I}$ is the unique dyadic
interval of length $2^{k_0}|I|$ which contains $I$.
We also observe that for different scales, the supports of
$\widehat{\widetilde{\Phi}^2_I}$ are disjoint.

Because of all these properties, we now observe that

\begin{equation}\label{32}
a^{(3)}_{J,k_0} =
\widetilde{B}_{k_0}(f_1, f_4), \Phi^3_J\rangle
\end{equation}
where

$$\widetilde{B}_{k_0}(f_1, f_4):= 2^{-k_0/2}
\sum_{I}
\frac{1}{|I|^{1/2}}
\langle f_1, \Phi_I^1\rangle
\langle f_4, \Phi_I^3\rangle
\widetilde{\Phi}_I^2.$$
As before, using now (\ref{32}) together with lemma \ref{l2} it follows that to prove
(\ref{30}) we just need to prove that

\begin{equation}\label{33}
\|\widetilde{B}_{k_0}(f_1, f_4)\widetilde{\chi}_{J_0}^{N'}\|_1
\lesssim
\|f_3 \widetilde{\chi}_{J_0}^N\|_p\cdot\|f_4 \widetilde{\chi}_{J_0}^N\|_q
\end{equation}
By scaling invariance, we may assume also as before that $|J_0| = 1$ and observe that then,
for every $I$ one has $|\widetilde{I}|\leq 1$.
Then, an argument similar to the one before allows us to reduce (\ref{33}) to

\begin{equation}\label{34}
\|\widetilde{B}_{k_0}(f_1, f_4)\chi_J\|_1
\lesssim
\|f_3 \widetilde{\chi}_{J}^N\|_p\cdot\|f_4 \widetilde{\chi}_{J}^N\|_q
\end{equation}
for every dyadic interval $J\subseteq \R$ of length $1$. It is thus sufficient to prove
(\ref{34}). We have, as before, several cases.

\underline{Case $II_{1a}$}: $\supp f_1$, $\supp f_4\subseteq 5J$.

Let $h\in\L^{\infty}$, $\|h\|_{\infty}\leq 1$. Then,

$$|\int_{\R} \widetilde{B}_{k_0}(f_1, f_4)(x) h(x) \chi_J(x) dx |\lesssim $$

$$2^{-k_0/2}\sum_{I}
\frac{1}{|I|^{1/2}}
|\langle f_1, \Phi_I^1\rangle|
|\langle f_4, \Phi_I^3\rangle|
|\langle h\chi_J,\widetilde{\Phi}_I^2\rangle|= $$

$$
\sum_I
|\langle f_1, \Phi_I^1\rangle|
|\langle f_4, \Phi_I^3\rangle|
\frac{|\langle h\chi_J,\widetilde{\Phi}_I^2\rangle|}{2^{k_0/2}|I|^{1/2}}
$$
and since now $\frac{\widetilde{\Phi}_I^2\rangle|}{2^{k_0/2}|I|^{1/2}}$ is 
$L^1$- normalized, the previous expression is smaller than

$$\sum_I|\langle f_1, \Phi_I^1\rangle|
|\langle f_4, \Phi_I^3\rangle| =
$$

$$\sum_I
\frac{|\langle f_1, \Phi_I^1\rangle|}{|I|^{1/2}}
\frac{|\langle f_4, \Phi_I^3\rangle|}{|I|^{1/2}}\cdot |I| =
$$

$$\int_{\R}
\sum_I
\frac{|\langle f_1, \Phi_I^1\rangle|}{|I|^{1/2}}
\frac{|\langle f_4, \Phi_I^3\rangle|}{|I|^{1/2}}
\chi_I(x) dx\lesssim
$$

$$\int_{\R}
\left(\sum_I\frac{|\langle f_1, \Phi_I^1\rangle|^2}{|I|}\chi_I(x)\right)^{1/2}\cdot
\left(\sum_I\frac{|\langle f_4, \Phi_I^3\rangle|^2}{|I|}\chi_I(x)\right)^{1/2} dx\lesssim
$$

$$\int_{\R} S(f_1)(x)\cdot S(f_4)(x) dx\lesssim
\|S(f_1)\|_p\cdot\|S(f_4)\|_q\lesssim
$$

$$\|f_1\|_p\cdot \|f_4\|_q \lesssim 
\|f_1\widetilde{\chi}^N_J\|_p\cdot
\|f_4\widetilde{\chi}^N_J\|_q
$$
using the fact that the square functions $S(f_1)$ and $S(f_4)$ are bounded on $L^r$ for
$1 < r <\infty$ and also the fact that we are in the Case $II_{1a}$.

\underline{Case $II_{1b}$}: Either $\supp f_1\subseteq (5J)^c$ or $\supp f_4\subseteq (5J)^c$

Assume as before that $\supp f_1\subseteq (5J)^c$. Then, decompose 
$\widetilde{B}_{k_0}(f_1, f_4)$ as

$$\widetilde{B}_{k_0}(f_1, f_4) = \widetilde{B}'_{k_0}(f_1, f_4) + 
\widetilde{B}''_{k_0}(f_1, f_4)$$
where

$$
\widetilde{B}'_{k_0}(f_1, f_4):=
2^{-k_0/2}
\sum_{\widetilde{I}\cap 5J\neq\emptyset}
\frac{1}{|I|^{1/2}}
\langle f_1, \Phi_I^1\rangle
\langle f_4, \Phi_I^3\rangle
\widetilde{\Phi}_I^2
$$
and

$$
\widetilde{B}''_{k_0}(f_1, f_4):=
2^{-k_0/2}
\sum_{\widetilde{I}\cap 5J = \emptyset}
\frac{1}{|I|^{1/2}}
\langle f_1, \Phi_I^1\rangle
\langle f_4, \Phi_I^3\rangle
\widetilde{\Phi}_I^2.
$$

If $h$ is as before, then we can write again

$$|\int_{\R}\widetilde{B}'_{k_0}(f_1, f_4)(x)  h(x) \chi_J(x) dx |\lesssim $$

\begin{equation}\label{35}
2^{-k_0}
\sum_{k=0}^{\infty}
\sum_{I; \widetilde{I}\cap 5J\neq\emptyset; |I| = 2^{-k}}
2^{2k}
|\langle f_1, 2^{-k/2}\Phi_I^1\rangle|
|\langle f_4, 2^{-k/2}\Phi_I^3\rangle|
|\langle h\chi_J, 2^{k_0/2}2^{-k/2}\widetilde{\Phi}_I^2\rangle|
\end{equation}
and we observe that the functions $2^{-k/2}\Phi_I^1$, $2^{-k/2}\Phi_I^3$ and
$2^{k_0/2}2^{-k/2}\widetilde{\Phi}_I^2$ are all $L^{\infty}$- normalized.
Then, we estimate (\ref{35}) by

$$
2^{-k_0}
\sum_{k=0}^{\infty}2^{2k}
\left(\sup_{I: \widetilde{I}\cap 5J\neq\emptyset; |I| = 2^{-k}}|
\langle f_1, 2^{-k/2}\Phi_I^1\rangle|\right)
\left(\sup_{I: \widetilde{I}\cap 5J\neq\emptyset; |I| = 2^{-k}}|
\langle f_4, 2^{-k/2}\Phi_I^3\rangle|\right)
\cdot$$

$$|\langle h\chi_J, \sum_{I: \widetilde{I}\in\widetilde{\I}; I\cap 5J\neq\emptyset; |I| = 2^{-k}}
\widetilde{\chi}_{\widetilde{I}}^N \rangle|.
$$
Since

$$|\langle h\chi_J, \sum_{I: \widetilde{I}\in\widetilde{\I}; I\cap 5J\neq\emptyset; |I| = 2^{-k}}
\widetilde{\chi}_I^N \rangle|\lesssim 2^{k_0},$$
the estimate follows as in the previous Case $I_{1b}$.

Finally, one can also write

$$|\int_{\R}\widetilde{B}''_{k_0}(f_1, f_4)(x)  h(x) \chi_J(x) dx |\lesssim $$

$$\sum_{I: \widetilde{I}\cap 5J = \emptyset}
\frac{1}{|I|^{1/2}}
|\langle f_1, \Phi_I^1\rangle|
|\langle f_4, \Phi_I^3\rangle|
|\langle h\chi_J,\widetilde{\Phi}_I^2\rangle|= $$

$$
\sum_{I: \widetilde{I}\cap 5J = \emptyset}
|\langle f_1, \Phi_I^1\rangle|
|\langle f_4, \Phi_I^3\rangle|
\frac{|\langle h\chi_J,\widetilde{\Phi}_I^2\rangle|}{2^{k_0/2} |I|^{1/2}}\lesssim
$$

$$\sum_{K: K\cap 5J = \emptyset}
\sum_{I: \widetilde{I} = K}
|\langle f_1, \Phi_I^1\rangle|
|\langle f_4, \Phi_I^3\rangle|
\frac{|\langle h\chi_J,\widetilde{\chi}^N_K\rangle |}{|K|}\lesssim
$$

$$
\sum_{K: K\cap 5J = \emptyset}
(\frac{\dist (K, J)}{|K|})^{-N'}
\sum_{I: \widetilde{I} = K}
|\langle f_1, \Phi_I^1\rangle|
|\langle f_4, \Phi_I^3\rangle|\lesssim
$$

$$\sum_{K: K\cap 5J = \emptyset}
(\frac{\dist (K, J)}{|K|})^{-N'}
\int_{\R}
\left(\sum_{I: \widetilde{I} = K}
\frac{|\langle f_1, \Phi_I^1\rangle|^2}{|I|}\chi_I(x)\right)^{1/2}\cdot
\left(\sum_{I: \widetilde{I} = K}
\frac{|\langle f_4, \Phi_I^3\rangle|^2}{|I|}\chi_I(x)\right)^{1/2} dx\lesssim
$$

$$\sum_{K: K\cap 5J = \emptyset}
(\frac{\dist (K, J)}{|K|})^{-N'}
\|f_1\widetilde{\chi}^N_K\|_p\cdot
\|f_4\widetilde{\chi}^N_K\|_q
$$
by using Lemma \ref{l2}. And this can be estimated further by

$$\sum_{K: K\cap 5J = \emptyset}
(\frac{\dist (K, J)}{|K|})^{-N'} (\dist (K, J)^{2N}
\|f_1\widetilde{\chi}^N_J\|_p\cdot
\|f_4\widetilde{\chi}^N_J\|_q.
$$
As in the Case $I_{1b}$ one observes that the sum

$$\sum_{K: K\cap 5J = \emptyset}
(\frac{\dist (K, J)}{|K|})^{-N'} (\dist (K, J)^{2N}
$$
is $O(1)$ if $N'$ is much bigger than $2N$ and so we obtain in the end the desired estimate.

\underline{Case $II_2$}: $j=3$.

This is actually easier, follows the same ideas and is left to the reader. This completes
the proof of Lemma \ref{l5}.

\section{Proof of Lemma \ref{l6}}

We are therefore left with proving Lemma \ref{l6} in order to complete the proof of our main
theorem.

\underline{Case $I$}: Estimates for $\energy^j_{3, \J}((a^{(3)})_J)$.

There are, as before, two subcases.

\underline{Case $I_1$}: $j\neq 3$.

Let $n\in\Z$ and $\bf{D}$ be so that the suppremum in Definition \ref{se1} is attained. Then,
since the intervals $J\in\bf{D}$ are all disjoint, we can write

$$\energy^j_{3, \J}((a^{(3)})_J)\sim 2^n (\sum_{J\in\bf{D}} |J| ) =
$$

$$= 2^n \|\sum_{J\in\bf{D}} \chi_J \|_1 = \|\sum_{J\in\bf{D}} 2^n \chi_J \|_{1, \infty}\lesssim
$$

\begin{equation}\label{36}
\left\|
\sum_{J\in\bf{D}}\frac{1}{|J|}
\|(\sum_{J'\subseteq J}\frac{|a^{(3)}_{J'}|^2}{|J'|}\chi_{J'}(x))^{1/2}\|_{1,\infty}
\chi_J \right\|_{1,\infty}.
\end{equation}
As in Section $9$, define the collection $\widetilde{\I}$ to be the set of all intervals
$I$ having the property that there exists $J\in\bf{D}$ and $J'\subseteq J$ with
$\omega^3_{J'}\cap \omega^2_I \neq \emptyset$ and $|\omega^3_{J'}| \leq |\omega^2_I|$. Then,
we observe as before that

\begin{equation}\label{36'}
a^{(3)}_{J'} = \langle B(f_1, f_4), \Phi^3_{J'}\rangle
\end{equation}
where $B(f_1, f_4)$ was defined by (\ref{24}). Using this fact together with Lemma \ref{l2}
one can majorize (\ref{36}) by

$$\left\|
\sum_{J\in\bf{D}}
\left(\frac{1}{|J|}\int_{\R} |B(f_1, f_4)| \widetilde{\chi}^N_J dx \right)
\chi_J \right\|_{1, \infty}\lesssim
$$

$$\|M ( B(f_1, f_4))\|_{1,\infty}\lesssim \|B(f_1, f_4)\|_1\lesssim
$$

\begin{equation}\label{37}
\sum_{I\in\tilde{\I}}
\frac{1}{|I|^{1/2}}
|\langle f_1, \Phi^1_I\rangle|
|\langle f_4, \Phi^3_I\rangle|
|\langle h, \Phi^2_I\rangle|
\end{equation}
for a certain $h\in L^{\infty}$, $\|h\|_{\infty}\leq 1$. Since $\frac{\Phi^2_I}{|I|^{1/2}}$
is an $L^1$- normalized function, it follows that (\ref{37}) is smaller than

\begin{equation}\label{38}
\sum_{I\in\widetilde{\I}}
|\langle f_1, \Phi^1_I\rangle|
|\langle f_4, \Phi^3_I\rangle|.
\end{equation}
Since both of the families $(\Phi^1_I)_I$ and $(\Phi^3_I)_I$ are {\it lacunary}, 
by a similar argument used to prove Proposition \ref{prop}, one can estimate the expression
(\ref{38}) by

$$\left(\size_{1,\widetilde{\I}}((\langle f_1, \Phi^1_I\rangle)_I)\right)^{1-\theta_1}
\left(\size_{2,\widetilde{\I}}((\langle f_1, \Phi^3_I\rangle)_I)\right)^{1-\theta_2}\cdot
$$

\begin{equation}\label{39}
\left(\energy_{1,\widetilde{\I}}((\langle f_1, \Phi^1_I\rangle)_I)\right)^{\theta_1}
\left(\energy_{2,\widetilde{\I}}((\langle f_1, \Phi^3_I\rangle)_I)\right)^{\theta_2}
\end{equation}
for any $0\leq \theta_1, \theta_2 < 1$ with $\theta_1 + \theta_2 = 1$ where
$\size_{1,\widetilde{\I}}((\langle f_1, \Phi^1_I\rangle)_I)$, 
$\size_{2,\widetilde{\I}}((\langle f_1, \Phi^3_I\rangle)_I)$,
$\energy_{1,\widetilde{\I}}((\langle f_1, \Phi^1_I\rangle)_I)$ and
$\energy_{2,\widetilde{\I}}((\langle f_1, \Phi^3_I\rangle)_I)$ are naturally defined as in
Definition \ref{se1}.

Using now the upper bounds for {\it sizes} and {\it energies} provided by Lemmas \ref{l3}
and \ref{l4}, (\ref{39}) can be estimated by

$$\left(\sup_I\int_{E_1}\widetilde{\chi}^N_I dx \right)^{1-\theta_1}
\left(\sup_I\int_{E_1}\widetilde{\chi}^N_I dx \right)^{1-\theta_2} |E_1|^{\theta_1}
|E_4|^{\theta_4}$$
which is the desired estimate.

\underline{Case $I_2$}: $j=3$.

This is easier. Pick again $n\in\Z$ and $\bf{D}$ so that the suppremum in Definition \ref{se1}
is attained. Then,

$$\energy^3_{3, \J}((a^{(3)})_J)\sim 2^n (\sum_{J\in\bf{D}} |J| ) =
$$

$$= 2^n \|\sum_{J\in\bf{D}} \chi_J \|_1 = \|\sum_{J\in\bf{D}} 2^n \chi_J \|_{1, \infty}\lesssim
$$

$$\left\|
\sum_{J\in\bf{D}}
\left(
\frac{1}{|J|}
\int_{\R}
|\sum_{I: \omega^3_J\cap\omega^2_I\neq\emptyset; |\omega^3_J| \leq |\omega^2_I|}
\frac{1}{|I|^{1/2}}
\langle f_1, \Phi^1_I\rangle
\langle f_4, \Phi^3_I\rangle
\Phi^2_I |
\widetilde{\chi}^N_J dx \right)
\chi_J
\right\|_{1,\infty}\lesssim
$$

$$\left\|
\sum_{J\in\bf{D}}
\left(
\frac{1}{|J|}
\int_{\R}
(\sum_{I}
|\langle f_1, \Phi^1_I\rangle|
|\langle f_4, \Phi^3_I\rangle|
\frac{\widetilde{\chi}_I^{N'}}{|I|})
\widetilde{\chi}^N_J dx \right)
\chi_J
\right\|_{1,\infty}\lesssim
$$

$$\left\|
M (\sum_{I}
|\langle f_1, \Phi^1_I\rangle|
|\langle f_4, \Phi^3_I\rangle|
\frac{\widetilde{\chi}_I^{N'}}{|I|})\right\|_{1,\infty}\lesssim
$$

$$\left\|
\sum_{I}
|\langle f_1, \Phi^1_I\rangle|
|\langle f_4, \Phi^3_I\rangle|
\frac{\widetilde{\chi}_I^{N'}}{|I|}\right\|_{1}\lesssim
$$

$$\sum_{I}
|\langle f_1, \Phi^1_I\rangle|
|\langle f_4, \Phi^3_I\rangle|
$$
and from here we can continue as before.

To obtain the estimates for $\energy^j_{3, k_0,\J}((a^{(3)})_J)$, one argues in the same way.
The $j=3$ case is identical to the corresponding previous one, while $j\neq 3$ follows also 
similarly. The only difference is that instead of (\ref{24}) one has to use (\ref{32}) and then
to observe that for every interval $I$, $\frac{\widetilde{\Phi}^2_I}{2^{k_0/2} |I|^{1/2}}$ 
is an $L^1$- normalized
function.

This ends our proof.

\end{document}